\documentclass[review,3p,10pt]{elsarticle}
\biboptions{sort&compress}

\usepackage{amssymb}
\usepackage{amsmath}
\usepackage{cases}
\usepackage{epstopdf}
\usepackage{fixltx2e}
\usepackage{mathtools}
\usepackage{array}
\usepackage{bm}
\usepackage{multirow}
\usepackage{makecell}
\usepackage{tabularx}
\usepackage{longtable}
\usepackage{float}
\usepackage{caption}
\usepackage{color}
\usepackage{algorithm}
\usepackage{algorithmic}
\usepackage{amsthm}
\usepackage{mathrsfs}
\usepackage{amssymb}

\usepackage{lineno}
\usepackage{booktabs}
\usepackage{caption}
\usepackage{graphicx}
\usepackage{float}
\usepackage{subcaption}
 \usepackage{pdflscape}
 \usepackage{adjustbox}
\usepackage[table,xcdraw]{xcolor}

\usepackage{url}

\newcommand{\tabincell}[2]{\renewcommand\arraystretch{0.8}\begin{tabular}{@{}#1@{}}#2\end{tabular}}

\newcommand{\variable}[2]{\textbf{#1} \> #2 \\}

\journal{Applied Energy}
\begin{document}
\begin{frontmatter}

%\title{Flexibility Assessment and Aggregation of Alkaline Electrolyzers under Dynamic Thermal and HTO Impurity Constraints for Renewable Power-to-Hydrogen Energy Management}
\title{Exploring the Optimal Size of Grid-forming Energy Storage in an Off-grid Renewable P2H System under Multi-timescale Energy Management}

\author[label1]{Jie~Zhu}
\author[label1]{Yiwei~Qiu\corref{cor1}}
\ead{ywqiu@scu.edu.cn}
\author[label1]{Yangjun~Zeng}
\author[label1]{Yi~Zhou}
\author[label1]{Shi~Chen}
\author[label1]{Tianlei~Zang}
\author[label1]{Buxiang~Zhou}
\author[label2]{Zhipeng~Yu}
\author[label2]{Jin~Lin}

\address[label1]{College of Electrical Engineering, Sichuan University, Chengdu, 610065, China}
\address[label2]{Department of Electrical Engineering, Tsinghua University, Beijing, 100087, China}
\cortext[cor1]{Corresponding author}

\begin{abstract}
  Utility-scale off-grid renewable power-to-hydrogen systems (OReP2HS), typically comprising photovoltaic plants, wind turbines, electrolyzers (ELs), and battery energy storage system (BESS), requires at least one component, generally the BESS, working with grid-forming ability to provide frequency and voltage references and regulate them through transient power support.
  However, existing designs of OReP2HS based on the energy management strategies (EMSs) with 5-min or hourly resolution fail to capture fast transients and may underestimate the BESS size required to ensure adequate grid-forming ability.
  This paper first proposes a framework of multi-timescale EMS covers from those for power system transient behaviors to second-level EL load adjustments and minute-level intra-day scheduling to coordinate renewable power, BESS, and ELs. Then, an iterative search procedure based on high-fidelity simulation (0.04 ms resolution) is employed to determine the cost-effective BESS size that satisfies grid-forming, long-term energy balancing over 8760 hours, and emergency support requirements.
  Case studies based on a planned OReP2HS project in Inner Mongolia, China, show that the proposed EMS yields a base-case LCOH of 33.212 CNY/kg (4.581 USD/kg), with CAPEX of BESS accounting for 17.83\% of total investment. The optimal BESS capacity represents 13.6\% of the rated hourly renewable output and shows a yearly degradation of 4.87\%.
  Sensitivity analysis reveals that reducing the electrolytic load adjustment time step from 90 to 5 seconds and increasing its ramping limit from 1\% to 10\% rated power per second, the BESS size decreases by 53.57\%, and the LCOH decreases to 25.458 CNY/kg (3.511 USD/kg). Considering the cost of designing and manufacturing utility-scale ELs with fast load regulation capability, a load adjustment time step of 5 to 10 seconds and a ramping limit of 4--6\% rated power per second are recommended to balance profitability and technological feasibility.
\end{abstract}


\end{frontmatter}

\section{Introduction}
\label{sec:intro}

\subsection{Background and Motivation}
\label{sec:Background and motivation}
{\color{black}In the global effort to reduce carbon emissions in electrical power, energy, transport, and chemical sectors, renewable power-to-hydrogen (ReP2H) technology has emerged as a promising pathway~\cite{li2021hydrogenSource}. According to the \textit{Global Hydrogen Review}~\cite{IEA2025Hydrogen}, the global installed capacity of electrolyzers (ELs) reached about 4.9 GW by September~2025 and is expected to exceed 26~GW by 2030 based on projects that have reached final investment decisions, with around 60\% of low-emission hydrogen projected to be produced by ReP2H.
%In the global effort to reduce carbon emissions, particularly in the electrical power, energy, transport, and chemical sectors, renewable power-to-hydrogen (ReP2H) technology has emerged as a promising solution. According to the \textit{China Hydrogen Development Strategy} \cite{meng2022ChinaHydrogenDevelopment}, the installed capacity of hydrogen electrolysers (ELs) in China reached 1.4 GW by 2022. Projections indicate this will increase to 100 GW in China and 35 GW in Europe by 2030, with 80\% of hydrogen demand met through ReP2H \cite{li2021hydrogenSource}.

Utility-scale ReP2H systems can be grid-connected or off-grid. Grid-connected systems can maintain steady hydrogen production by exchanging power with the utility grid to manage renewable power fluctuations. However, ensuring all hydrogen produced meets green certification is challenging \cite{dufo2024optimisation,he2021transformation}, and remote areas with abundant renewable resources, such as islands and deserts, often lack grid access. In contrast, off-grid ReP2H systems (OReP2HSs) offer flexible design according to the spatial distribution of renewable resources and avoid grid-exchange constraints. Supported by policy initiatives, numerous demonstration projects are now under construction~\cite{zeng2025planning,zeng2024investment,ChinaP2M}.}
%Conversely, off-grid ReP2H systems (OReP2HS) can be designed flexibly based on the spatial distribution of the renewable energy resources. The stringent constraints and regulations on power exchange with the utility grid are also avoided. Benefiting from these advantages, OReP2HS has gained policy support, leading to many demonstration projects under construction \cite{zeng2025planning,yu2025novel,ChinaP2M}.

The OReP2HS typically comprises photovoltaic (PV) plants, wind turbines (WTs), ELs, and a battery energy storage system (BESS) \cite{abdin2015solar}. Alternating current (AC) paradigms are often preferred due to the immature market for direct current (DC) grid equipment \cite{he2020review,gallardo2022assessing}. In the off-grid mode, at least one component in the OReP2HS must operate under grid-forming control (such as voltage/frequency droop and virtual synchronous generator control) to provide frequency and voltage references and regulate them through transient power support and short-term energy balance regulation. However, the development of grid-forming capabilities for WTs, PVs, and ELs is still in the experimental stage \cite{pawar2021gridFormingControl,tavakoli2023gridFormingControlALE,qi2021gridFormingControlWT}. Using BESS as the grid-forming component appears to be the most feasible solution \cite{zhao2021controlGFMBESS,zuo2021performanceGFMBESS}, but poses challenges in determining its configuration.

First, the grid-forming capability of OReP2HSs inherently involves fast transient dynamics that cannot be captured by existing low-resolution models using 5-minute or hourly time steps~\cite{dysko2020testingCharateristicGFM,gerini2022optimalGFM}. {\color{black}Although BESS sizing based on such low-resolution models can achieve power balance at discrete time points, it may fail to satisfy transient stability requirements, since the intra-interval dynamics of frequency and voltage are not represented. This omission can lead to an underestimation of the BESS capacity required to provide reliable grid-forming support, thereby compromising system operational stability. To address this, a high-resolution simulation model and an energy management system (EMS) capable of capturing and regulating transient behaviors are essential for determining the minimum BESS size that ensures frequency and voltage stability while minimizing investment costs.}

Second, ELs offer a wide operating range and near-second-level response capability \cite{qiu2023extended}. An appropriately designed EMS that coordinates PV, WTs, BESS, and ELs to achieve second-level source-load power balancing can significantly reduce the required BESS capacity, thereby enhancing the overall economic viability of OReP2HS. However, such an EMS must span multiple timescales from grid-forming control and second-level EL load adjustment to intra-day/day-ahead scheduling, that are lacking in existing studies.

To address these limitations, this research proposes a comprehensive multi-timescale EMS that covers from transient grid-forming behaviors to rolling intra-day scheduling of ELs. Leveraging high-fidelity simulation models with 0.04 ms time resolution, the BESS size is optimized to minimize the levelized cost of hydrogen (LCOH) while ensuring grid-forming ability, year-long energy balance, and emergency support capability. Besides, the influence of various factors on BESS sizing is analyzed. A review of relevant literature is presented in Section \ref{sec:review}, and the contributions of this study are summarized in Section \ref{sec:contribution}.

\subsection{Literature Review}
\label{sec:review}
The EMS in the existing works for supporting the production simulation to optimize the configuration of an OReP2HS can be divided into three categories: rule-based, {\color{black}mathematical} optimization-based, or directly using EMS integrated in commercial planning softwares.

Most existing studies design rule-based EMS for OReP2HS operation, a typical logic of such EMS is as follows: when the available renewable power exceeds the minimum load of the ELs, the ELs are started, and their load follows the renewable power with a 5-min or hourly fixed time step. Once the ELs reach their maximum capacity, any surplus power is used to charge the BESS. Conversely, if the combined power from renewable power and the BESS becomes insufficient to sustain the ELs at their minimum load, the ELs are shut down. For example, Marocco et al. \cite{marocco2020study}  designed a rule-based EMS for an OReP2HS to maintain energy autonomy in remote microgrids. When renewable power exceeds load demand, the BESS is charged first. Once the BESS reaches its maximum state of charge (SOC), ELs are started up, and surplus power is directed to hydrogen production. Based on the same EMS, Gandiglio et al. \cite{gandiglio2022life} assessed the environmental benefits of powering an island village with an OReP2HS. Abdulrahman et al. \cite{al2022hydrogen} used a similar OReP2HS to produce hydrogen for refueling hydrogen-based vehicles in an off-grid residential area in Dhahran, Saudi Arabia, and obtained a minimal LCOH of 36.32 USD/kg. Under a rule-based EMS, metaheuristic algorithms are usually employed to optimize the sizes of system components. For instance, Lu et al. and Marocco et al. \cite{lu2023modeling, marocco2022role} applied particle swarm optimization (PSO) to optimize the capacities of PVs, WTs, BESS, ELs, and hydrogen tanks (HTs) in an OReP2HS. Charged system search algorithms \cite{kohole2024optimization}, genetic and gray wolf optimization \cite{mohseni2020economic}, and NSGA-\uppercase\expandafter{\romannumeral2} \cite{xu2020data} are also used for sizing the OReP2HS. However, rule-based EMS generally follows static rules, which can not addressing uncertainties of renewable power. Besides, it covers only one time-scale, often hourly, which can not capture the transient dynamic of grid-forming BESS and the near-second-level load response capability of ELs.

%metaheuristic algorithms are usually employed to optimize the sizes of the system components. For example, \cite{marocco2020study,marocco2022role,marocco2021optimal} designed a rule-based EMS for an OReP2HS to maintain energy autonomy in remote microgrids. When renewable power exceeds load demand, the BESS is charged first. Once the BESS reaches its maximum state of charge (SOC), ELs are started up, and surplus power is directed to hydrogen production. The capacities of PVs, WTs, BESS, ELs, and hydrogen tanks (HTs) are optimized using particle swarm optimization (PSO) with hourly time resolution \cite{marocco2022role}. A techno-economic analysis, considering the degradation cost of batteries and fuel cells, is presented in \cite{marocco2021optimal}. Gandiglio et al. \cite{gandiglio2022life} assessed the environmental benefits of powering an island village with an OReP2HS using the EMS proposed in \cite{marocco2020study}.
%
%The OReP2HS also plays a crucial role in achieving energy autonomy for green buildings \cite{lu2023modeling}, standalone communities \cite{kohole2024optimization,mohseni2020economic}, and off-grid industrial parks \cite{xu2020data}. The system sizes are optimized using PSO \cite{lu2023modeling}, charged system search algorithms \cite{kohole2024optimization}, genetic algorithms \cite{mohseni2020economic}, and gray wolf optimization \cite{mohseni2020economic}, as well as NSGA-\uppercase\expandafter{\romannumeral2} \cite{xu2020data} based on rule-based EMSs.

Other researchers designed mathematical optimization-based EMSs and determined component sizes for the OReP2HS through deterministic optimization. In these studies, OReP2HS operates following the results of a day-ahead/intra-day optimal scheduling. For instance, Oyewole et al. \cite{oyewole2024optimal} proposed {\color{black}an integrated} mixed-integer linear programming (MILP) approach for designing OReP2HS in insular communities. Linear programming (LP) was utilized for planning an OReP2HS to supply a research facility in isolated regions \cite{viole2023renewable}. For transportation and industrial applications, Wang et al. \cite{wang2023optimising} focused on green hydrogen production for green ammonia synthesis, employing MILP for simulations and optimizations. Ibagon et al. \cite{ibagon2023techno} investigated the optimal size of OReP2HS for export-oriented green hydrogen production in Uruguay, using sequential quadratic programming (SQP) and predicting a decrease in LCOH from 3.5 USD/kg to 2.3 USD/kg by 2030 with technological maturity. Yang et al. \cite{yang2020planning} designed an off-grid hydrogen supply chain in Fujian, China, using chance-constrained programming to address planning and operation problems. Shao et al. \cite{shao2023risk} proposed bi-level stochastic programming for dispatch and size determination of OReP2HS, while Pang et al. \cite{pang2022integrated} suggested a similar bi-layer framework using mixed integer quadratic constrained programming (MIQCP) for scheduling.

Studies adopt commercial planning software to design OReP2HS is similar to mathematical optimization-based methods, since the EMSs integrated in softwares are often MILP-based model. Babaei et al. \cite{babaei2022optimization} used OReP2HS to supply power to energy-stressed islands in Eastern Canada, employing HOMER to find the most cost-effective configuration. Abdin et al. \cite{abdin2024feasibility} used HOMER PRO to optimize the sizes of PV, WTs, ELs, HTs, and BESS. Similarly, HOMER/HOMER PRO was used to find cost-effective configurations of OReP2HS in remote communities \cite{babatunde2022off} and villages \cite{al2017techno}. DECAPLA was applied for simulation in the design of OReP2HS for independent communities, considering an MILP-based EMS \cite{nastasi2023renewable}.

However, EMS in the above optimization- or software-based studies only considered hourly time resolution, overlooking renewable power fluctuations on shorter timescales. To this end, Ib{\'a}{\~n}ez-Rioja et al. \cite{ibanez2022simulation,ibanez2023off} used 300-second-resolution data from a Finnish wind farm to improve the time resolution for OReP2HS production simulations and optimal configuration, obtaining a more precise LCOH of 2 Euro/kg by 2030 in southeastern Finland. Zheng et al. \cite{ZHENG2023113763} used 5-minute-resolution data to evaluate the cost-effective component size for an off-grid wind power-to-hydrogen system, finding an LCOH range of 1.66--3.61 USD/kg in wind-rich areas in China and Denmark. Nevertheless, EMS with 5-min resolution still can not meet the requirement of simulation covers a transient timescale for the OReP2HS designing.

\begin{figure}[t]
	\centering
	\includegraphics[width=5.5in]{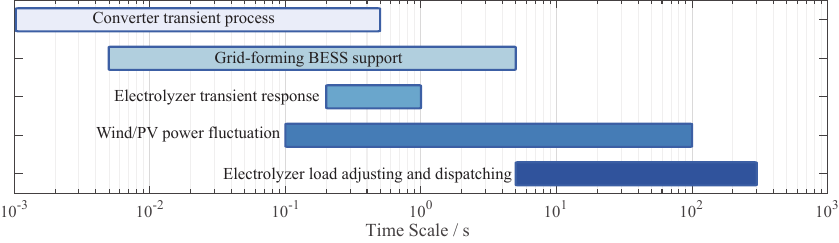}
			\vspace{-6pt}
	\caption{Physical processes and time scale distribution of a typical OReP2HS.}
	\label{fig:TimeScale}
\end{figure}

The reviewed studies related to the planning and designing of OReP2HS are summarized in Table \ref{tab:literature}. These studies mainly focus on the 8760-hour or day-ahead/intraday energy balance constraints and scheduling of the OReP2HS using 5-min or hourly time resolution models, which only consider power balance at discrete time steps, failing to capture the dynamic interactions among PV, WTs, BESS, and ELs at shorter time scales.
{\color{black}This results in the following shortcomings:

First, beyond the minute- or hour-level source-load power and energy balance, the OReP2HS, as an independent power system, also exhibits the fundamental interactions dynamic of generation, grid, and loads~\cite{qiu2024technological}. These include the electrical energy conversion and transient control processes of converters of WT/PV/EL/BESS ~\cite{9265486}, as well as the frequency and voltage dynamic of the system~\cite{lu2025stability}, which typically occur on millisecond to second time scales. Besides, the transient electrochemical response at the electrode-electrolyte interface of EL making its load response exhibit a first one order dynamic response with a time constant of approximately 0.2 to 1 second~\cite{cheng2025power,sha2023low}. Fig.~\ref{fig:TimeScale} gives the overall timescales involved a typical OReP2HS, ranging from milliseconds to hours.  Therefore, existing studies design the OReP2HS employ 5-min or hourly time resolution models cannot accurately reflect its real dynamic behavior, which may compromise the feasibility of the resulting design outcomes in real-word applications.}
%This results in two shortcomings in : 1) underutilization of the ELs’ fast load regulation capability, and 2) inadequate consideration of the grid-forming requirements of the BESS.

\begin{table}[tb]
  \scriptsize
  \renewcommand{\arraystretch}{1.1}
  \caption{Summary of the literatures related to the planning and designing of OReP2HS.}
  \label{tab:literature}
  \vspace{-6pt}
  \centering
  \begin{adjustbox}{width=\textwidth}
  \begin{tabular}{ccccccc}
    \hline \hline
    \makebox[0.08\textwidth][c]{\makecell[c]{Literature}} &
    \makebox[0.08\textwidth][c]{\makecell[c]{Components}} &
    \makebox[0.07\textwidth][c]{\makecell[c]{8760-hour\\energy balance}} &
    \makebox[0.12\textwidth][c]{\makecell[c]{Day-ahead/Intra-\\day scheduling}} &
    \makebox[0.09\textwidth][c]{\makecell[c]{Frequency/\\voltage security}} &
    \makebox[0.10\textwidth][c]{\makecell[c]{Resolution\\of simulation}} &
    \makebox[0.12\textwidth][c]{\makecell[c]{Method/\\Software}} \\
    \hline
    \makecell[c]{\cite{marocco2020study}, \cite{gandiglio2022life}} &
    \makebox[0.1\textwidth][c]{\makecell[c]{PV/WT/BT$^1$/\\EL/HT}} &
    \checkmark & & & Hourly & \makecell[c]{Iterative searching} \\

    \cite{al2022hydrogen} &
    \makebox[0.1\textwidth][c]{\makecell[c]{PV/WT/BT/EL}} &
    \checkmark & & & Hourly & \makecell[c]{Iterative searching} \\

    \makecell[c]{\cite{lu2023modeling}, \cite{marocco2022role}} &
    \makebox[0.1\textwidth][c]{\makecell[c]{PV/WT/BT/EL/HT}} &
    \checkmark & & & Hourly & PSO \\

    \cite{kohole2024optimization} &
    \makebox[0.1\textwidth][c]{\makecell[c]{PV/WT/EL/DG$^2$}} &
    \checkmark & & & Hourly & \makecell[c]{Metaheuristic} \\

    \cite{mohseni2020economic} &
    \makebox[0.1\textwidth][c]{\makecell[c]{PV/WT/SC$^3$/EL/HT/\\Micro-hydro}} &
    \checkmark & & & \makecell[c]{Hourly and\\monthly} & \makecell[c]{Metaheuristic} \\

    \cite{xu2020data} &
    PV/WT/EL/HT &
    \checkmark & & & Hourly & \makecell[c]{NSGA-II} \\

    \cite{oyewole2024optimal} &
    \makebox[0.1\textwidth][c]{\makecell[c]{PV/WT/EL/HT}} &
    \checkmark & \checkmark & & Hourly & MILP \\

    \cite{viole2023renewable} &
    \makebox[0.1\textwidth][c]{\makecell[c]{PV/WT/BT/EL/HT/DG}} &
    \checkmark & \checkmark & & Hourly & LP \\

    \cite{wang2023optimising} &
    \makebox[0.1\textwidth][c]{\makecell[c]{PV/WT/BT/EL}} &
    \checkmark & \checkmark & & Hourly & MILP \\

    \cite{ibagon2023techno} &
    \makebox[0.1\textwidth][c]{\makecell[c]{PV/WT/BT/EL}} &
    \checkmark & \checkmark & & Hourly & SQP \\

    \cite{yang2020planning} &
    \makebox[0.1\textwidth][c]{\makecell[c]{PV/WT/BT/EL}} &
    \checkmark & \checkmark & & Hourly & \makecell[c]{Chance constrained\\programming} \\

    \cite{shao2023risk} &
    \makebox[0.1\textwidth][c]{\makecell[c]{PV/WT/BT/EL/HT}} &
    \checkmark & \checkmark & & Hourly & \makecell[c]{Stochastic\\programming} \\

    \cite{pang2022integrated} &
    \makebox[0.1\textwidth][c]{\makecell[c]{PV/WT/BT/EL/Dispenser}} &
    \checkmark & \checkmark & & Hourly & MIQCP \\

    \makecell[c]{\cite{babaei2022optimization}, \cite{babatunde2022off}, \cite{al2017techno}} &
    \makebox[0.1\textwidth][c]{\makecell[c]{PV/WT/BT/EL/HT}} &
    \checkmark & \checkmark & & Hourly & HOMER \\

    \cite{abdin2024feasibility} &
    \makebox[0.1\textwidth][c]{\makecell[c]{PV/WT/BT/EL}} &
    \checkmark & \checkmark & & Hourly & HOMER Pro \\

    \cite{nastasi2023renewable} &
    \makebox[0.1\textwidth][c]{\makecell[c]{PV/WT/BT/EL/HT}} &
    \checkmark & \checkmark & & Hourly & DECAPLA \\

    \cite{ibanez2022simulation} &
    \makebox[0.1\textwidth][c]{\makecell[c]{PV/WT/BT/EL}} &
    \checkmark & & & 5-min & PSO \\

    \cite{ibanez2023off} &
    \makebox[0.1\textwidth][c]{\makecell[c]{PV/WT/BT/EL}} &
    \checkmark & & & 5-min & \makecell[c]{Metaheuristic} \\

    \cite{ZHENG2023113763} &
    \makebox[0.1\textwidth][c]{\makecell[c]{WT/EL}} &
    \checkmark & & & 5-min & \makecell[c]{Iterative searching} \\
    \hline
  \bf This work                                     & \bf PV/WT/BT/EL       & \bf\checkmark  & \bf\checkmark     & \bf \checkmark     & \makebox[0.1\textwidth][c]{\tabincell{c}{\bf 0.04 ms}}  &  \makebox[0.1\textwidth][c]{\tabincell{c}{ \bf High-fidelity \\ \bf  simulation-based \\ \bf searching}}               \\
  \hline \hline
  \multicolumn{7}{l}{\makecell[l]{$^1$Battery; $^2$Diesel generator; $^3$Super capacitor. }}   \\

  \end{tabular}
  \end{adjustbox}
\end{table}

Second, many researchers have investigated the use of ELs for frequency regulation in the power grid \cite{cheng2022coordinated,dozein2021fast,dozein2022virtual}, demonstrating the capability of the ELs for load regulation within seconds. This capability can smooth renewable power fluctuations, reducing BESS capacity requirements. However, this second-level capability is not considered in existing works with 5-min or hourly resolution, leading to a potential overestimation of the optimal BESS size.

Third, since the absence of grid energy support, grid-forming BESS provides necessary frequency and voltage references for an OReP2HS and adjusts its output power instantaneously to maintain them within thresholds \cite{rosso2021gridforming}. If the power imbalance exceeds the BESS capacity, the OReP2HS may destabilize immediately or shut down due to the loss of operating references \cite{lasseter2019gridforming}. This grid-forming dynamic of BESS occurs on the time scale of milliseconds and involves complex interactions among the control loops of PV, WTs, and ELs. {\color{black} Simulations reported in the literature reviewed in Table \ref{tab:literature} with minute-level to hourly resolution, cannot truly reveal the real grid-forming process of an OReP2HS and consider frequency and voltage security constraints in designing, risking the underestimation of the BESS size required to provide reliable grid-forming support and erroneous estimation of LCOH.}
%Thus, the LCOH ranges of 1.6--3.5 USD/kg reported in the literature \cite{ibagon2023techno,ibanez2022simulation,ibanez2023off,ZHENG2023113763} appear to be overly optimistic.

%In summary,

\subsection{Contributions}
\label{sec:contribution}
To fill the aforementioned gaps, {\color{black}this paper proposes high-fidelity simulation-based grid-forming BESS optimal sizing method, which employee an enhanced optimization-based EMS with a comprehensive high-fidelity model for the design of OReP2HS while considering frequency and voltage security constraints.} Unlike existing studies, the proposed EMS operates at a significantly higher time resolution, enabling the capture of system transient behaviors. It integrates multi-layer control and scheduling across various timescales to ensure both economic and operational security. A planned utility-scale OReP2HS project in Inner Mongolia, China, is used as a case study to provide practical insights for industrial implementation. The main contributions are summarized as follows:
\begin{enumerate}
	\item	We first propose a multi-timescale EMS for OReP2HS covering a range of millisecond-level transient power balance to 4-hour intra-day scheduling. It coordinates MILP-based rolling scheduling, second-level EL load-following control, an emergency handling strategy to address contingencies such as WT/PV/EL tripping and grid-forming control of BESS.
	
	 \item	A comprehensive model is developed to support the high-fidelity production simulation of an OReP2HS with 0.04-millisecond time resolution. {\color{black}This model integrates the proposed EMS with electromagnetic transient simulation models, which, for the first time, captures the millisecond-level transient power balance among PV, WTs, ELs, and BESS, and frequency and voltage dynamic during the BESS sizing.}

	 \item  The optimal size of the grid-forming BESS in an OReP2HS is evaluated through iterative searching to achieve minimal LCOH. {\color{black}Constraints for long-term (8760 hours) energy balance, grid-forming ability (frequency and voltage securities), and continuous operation under emergency conditions are simultaneously considered for the first time.}
	
	\item  The impact of the proposed EMS on the optimal BESS size and LCOH is assessed. The findings indicate that the base-case optimal BESS capacity is 6.8 MW/3.4 MWh, yielding an LCOH of 33.212 CNY/kg (4.581 USD/kg). {\color{black}Further sensitivity-based techno-economic analysis reveals that both the required BESS size and the LCOH can be further reduced by enhancing the EL ramping capability and decreasing the load adjustment time step.}
\end{enumerate}

The rest of this paper is organized as follows. The configuration and component model of the OReP2HS are presented in Section \ref{sec:System description}. In Section \ref{sec:EMS}, the multi-timescale EMS is described in detail. Section \ref{sec:BatOptim} presents the method for evaluating the optimal BESS size. The results and sensitivity analysis are given in Section \ref{sec:CaseStudy}. Finally, Section \ref{sec:Conclusion} presents the conclusion and future prospects.

\section{System Description}
\label{sec:System description}

\subsection{ Configuration of the OReP2HS}
\label{sec:Config}

\begin{figure}[t]
	\centering
	\includegraphics[width=5.5in]{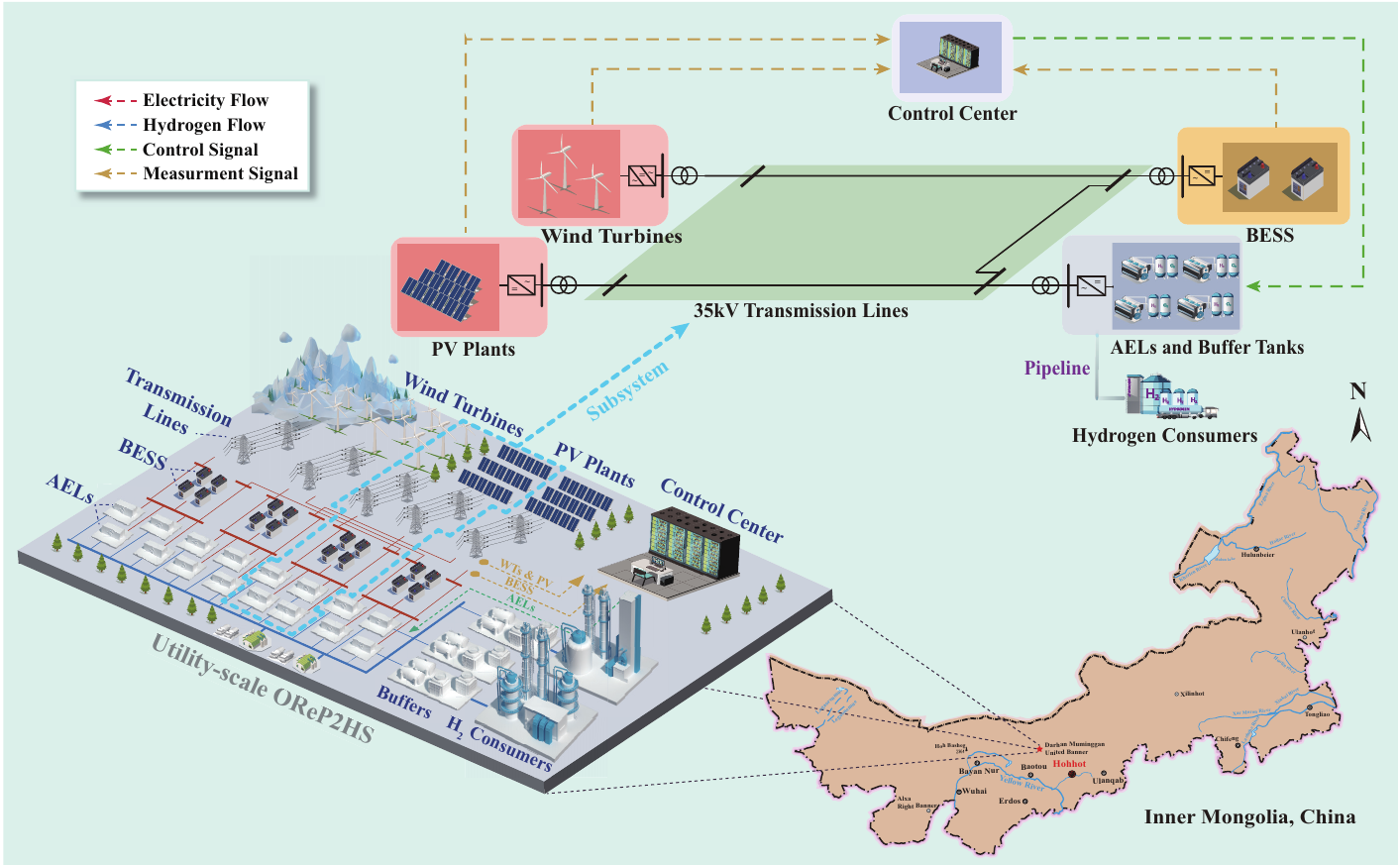}
			\vspace{-6pt}
	\caption{Schematic diagram of the studied OReP2HS located in Inner Mongolia, China.}
	\label{fig:config}
\end{figure}

The OReP2HS studied in our research is based on a realistic project planned in Inner Mongolia, China. Its schematic diagram is shown in Fig. \ref{fig:config}. This utility-scale OReP2HS comprises multiple subsystems that are electrically independent. The produced hydrogen is collected by pipelines and fed to industrial consumers in a chemical park.

{\color{black}Since the subsystems operate electrically independently without transient interactions, this paper takes one single subsystem as the research object.} It is planned to install three 6.25 MW WTs from {\color{black}\emph{Goldwind Sci \& Tech}~\cite{goldwind_equipment}} and a 6.25 MWp PV plant from {\color{black}\emph{LONGi Green Energy Technology}~\cite{longi_modules}}. The hydrogen production plant comprises four alkaline ELs (AELs) from {\color{black}\emph{Peric Hydrogen Technologies}~\cite{peric2025}}. Each AEL has a rated power of 5 MW and a maximum overload of 20\% of the rating. A BESS from {\color{black}\emph{SUNGROW}~\cite{noauthor_product_nodate}} is taken as the grid-forming equipment, and the evaluation of its size, characterized by the battery capacity and charging/discharging rate, is the goal of this work. All devices are connected by 35 kV AC transmission lines and operate at a rated frequency of 50 Hz.

\subsection{Requirements for Sizing the Grid-forming BESS}
\label{sec:Requorement}

The size of the grid-forming BESS need to ensure economic viability of the OReP2HS while meet the following requirements:
\begin{enumerate}
 \item{\bf Grid-forming ability:} The BESS must provide frequency and voltage references and offer transient power support to maintain these parameters within limits. This is a prerequisite for the stable operation of an OReP2HS. While increasing battery capacity and charging/discharging rates enhances grid-forming ability, it also raises investment, operation, and maintenance costs.

 \item{\bf Continuous operation during emergencies:} The OReP2HS must immediately balance power shortages or surpluses to maintain frequency and voltage security during emergencies, such as the tripping of a WT, PV, or EL. This capability is essential for preventing unnecessary blackouts and ensure a reliable hydrogen supply.

 \item{\bf Long-term energy balance:} The OReP2HS must coordinate PV, WTs, BESS, and ELs to maintain year-long (8760 h) energy balance while avoiding the BESS being fully charged/discharged to maintain the power regulation ability. Note that under an effective load scheduling strategy, the regulation range of electrolytic load is far beyond that of the BESS. A full utilization of this capability can reduce the size of the BESS and the LCOH.
\end{enumerate}

The grid-forming ability and continuous operation during emergencies involve fault (e.g., $N$-  1 contingency) analysis and complex dynamic interactions on a millisecond timescale between PV, WTs, BESS, and ELs. A refined simulation model is necessary to test them accurately. On the other hand, when testing long-term energy balance, a trade-off between time efficiency and accuracy is needed, making simplified models more suitable. Therefore, we develop both refined and simplified simulation models in Section \ref{sec:Model}. Then, a multi-timescale EMS of the OReP2HS is presented in Section \ref{sec:EMS} to support the simulations.

\subsection{Simulation Models of the OReP2HS}
\label{sec:Model}

Refined and simplified simulation models are developed to test the three requirements listed in Section \ref{sec:Requorement}. These models are derived from real-world devices as introduced in Section \ref{sec:Config} and are described below. Detailed topologies, controllers, and parameters are presented in \ref{sec:appA}.

The refined model is developed in the electromagnetic transient form and established in \emph{Matlab/Simulink}. For WTs and the PV plant, maximum power point tracking (MPPT) control is integrated into the control loops. The influence of the ambient temperature on the PV plant output is also considered. For the AELs, detailed models considering electrolysis stacks and rectifiers are developed. The rectifier includes vector-orientation-controlled IGBT-PWM-based AC/DC and Buck/Boost-based DC/DC converters.

{\color{black}The stack is modeled as an equivalent circuit consisting of a nonlinear resistance in series with a counter electromotive force to reproduce its U-I feature, which capture the quasi-steady-state electrochemical behavior \cite{iribarren2023dynamic}. In addition, a load ramping limit mechanism is introduced to emulate the impact of the electrical double-layer (EDL) on stack load adjustment, thereby representing its transient electrochemical response \cite{cheng2025power}.} The grid-forming BESS adopts voltage/frequency (V/F) droop control. A generic model represents the battery as a controllable ideal DC source with an equivalent resistance \cite{farrokhabadi2017battery}.

The simplified model has the same topologies and parameters as the refined model but differs with regard to power electronic converters modeling. In the refined model, converters are modeled with detailed circuit topology and switching dynamics of the IGBT/diode pairs controlled by PWM signals. For the simplified model, the switching dynamics within the converters are disregarded. Instead, a switching-function model directly controlled by the reference voltage signals \cite{MathWorksTwoLevelConverter} is used to improve efficiency in testing long-term energy balance.

The hybrid modeling combining both high-fidelity and simplified models could better simulate real-world operating conditions and addresses the limitations of low-resolution (e.g., with 5-min time steps) models.

\section{Multi-timescale Energy Management System}
\label{sec:EMS}

\subsection{Framework}
\label{sec:Framework}

\begin{figure}[tb]
  \centering
  \includegraphics[width=5.8in]{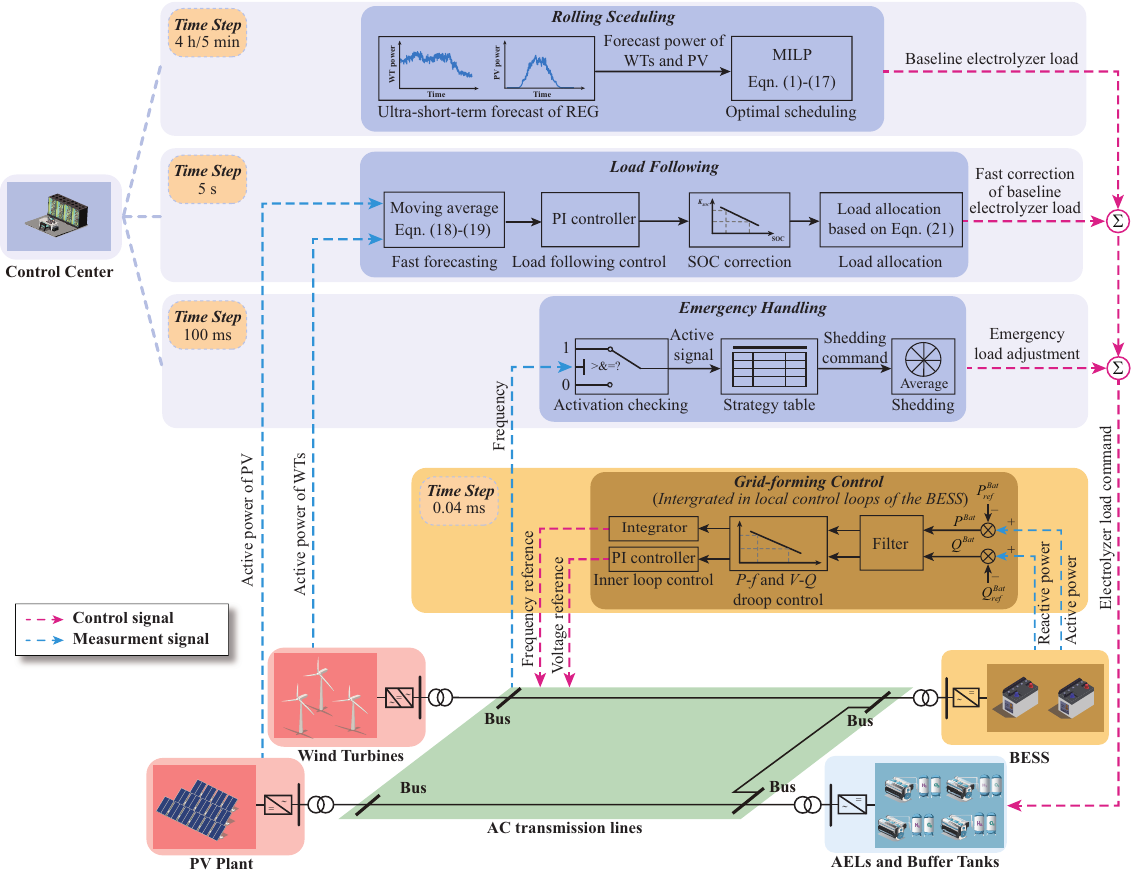}
  		\vspace{-3pt}
  \caption{Framework of the proposed multi-timescale EMS.}
  \label{fig:EMS}
\end{figure}

\begin{figure}[tb]
  \centering
  \includegraphics[width=5.8in]{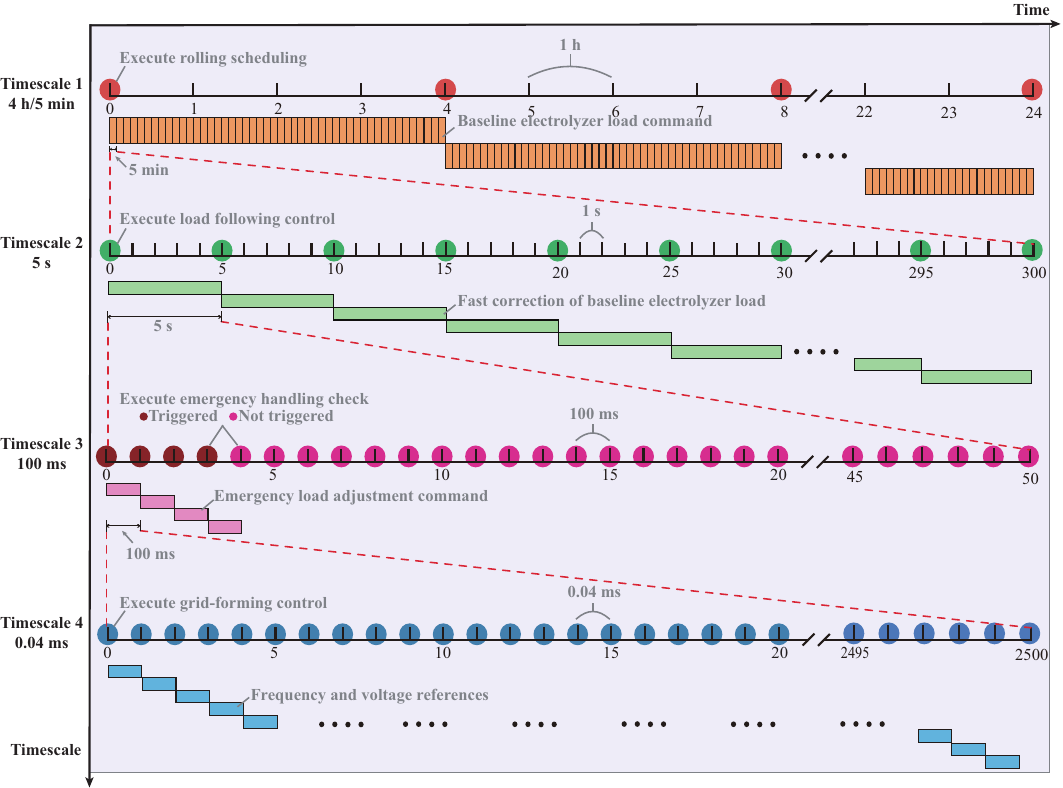}
  		\vspace{-3pt}
  \caption{{\color{black}Temporal coordination among the four sub-strategies of the proposed multi-timescale EMS.}}
  \label{fig:TimeLogic}
\end{figure}

The framework of the proposed multi-timescale EMS is illustrated in Fig. \ref{fig:EMS}. WTs and the PV plant operate with the built-in MPPT controllers. The grid-forming BESS operates under V/F droop control to provide voltage and frequency references to the OReP2HS. The control center regulates the electrolytic load to track the renewable power.

%The timescale of the proposed EMS covers the range from 0.04 ms to 4 h, and four sub-strategies are included. The strategies are coordinated to achieve  stable and economic operation of the OReP2HS.
%\begin{enumerate}
%  \item{\bf Rolling Scheduling (from 5 min to 4 h):} The rolling scheduling (RS) aligns the hydrogen production plan with the renewable generation profile based on forecasts. It utilizes an MILP to determine the start-up/shut-down/standby states and baseline load for each EL in every upcoming 4 h. Its time step is typically from 5 to 15 min.
%
%  \item{\bf Load Following (typically 5 s):} The load following (LF) is similar to automatic generation control (AGC) in power systems. It adjusts the load of AELs every 5 s. The goal is to ensure that the electrolytic load closely tracks the renewable power, thereby reducing the regulation demand on the grid-forming BESS.
%
%  \item{\bf Emergency Handling (100 milliseconds):} The emergency handling (EH) ensures the continuous operation of the OReP2HS during emergencies. It monitors system frequency and power balance and adjusts the power of PV, WTs, BESS, and ELs in the event of contingencies, such as WT/PV/EL tripping.
%
%  \item{\bf Grid-forming Control (0.04 milliseconds):} Under grid-forming control, the BESS provides frequency and voltage references in the electrical network, offering transient power support and energy balance regulation to maintain stability.
%\end{enumerate}

Four sub-strategies are included in the proposed EMS, and Fig. \ref{fig:TimeLogic} illustrates the temporal coordination among them. These strategies are coordinated to achieve  stable and economic operation of the OReP2HS.
\begin{enumerate}
  \item{\bf Rolling Scheduling (4 h/5 min):} The rolling scheduling (RS) operates at the top level, executed every 4 hours. At each interval, it takes the power forecasts provided by renewable power plants, current load power and start-up/shut down/stand by states of ELs as the initial conditions to optimizes the unit commitment and baseline power trajectory of each EL for the next 4 hours, with a scheduling time step of 5 minutes.

  \item{\bf Load Following (typically 5 s):} Within each 5-min scheduling interval, a load-following (LF) control is executed every 5 seconds, providing fast corrections to the baseline power setpoint of each EL. It is similar to automatic generation control (AGC) in power systems, to ensure that the electrolytic load closely tracks the renewable power, thereby reducing the regulation demand on the grid-forming BESS.

  \item{\bf Emergency Handling (100 milliseconds):}  Within each 5-second load-following cycle, an emergency handling (EH) check is conducted every 100 milliseconds. If triggered, a supplementary power shedding command is issued to adjust the EL load accordingly, ensuring power balance and frequency stability of the OReP2HS during emergencies.

  \item{\bf Grid-forming Control (0.04 milliseconds):} V/F droop control for grid-forming operates as the foundational control mechanism, running every 0.04 milliseconds within each emergency check cycle. Under grid-forming control, the BESS provides voltage and frequency references to support the OReP2H system operation and offers transient power support and energy balance regulation to keep them stable.
\end{enumerate}

\subsection{Rolling Scheduling}
\label{sec:RS}
The RS module determines the operational states (start-up, shut-down, or standby) and baseline load of each EL over a 4-hour horizon, based on ultra-short-term renewable power forecasts provided by the power plants. This module uses a MILP approach with a time step $\Delta t_{\mathrm{RS}}=5$ minutes. The objective is to maximize hydrogen production revenue, as described by (\ref{eq-1}),
\begin{equation}\label{eq-1}
  \min_{\boldsymbol{x}}\ \sum\nolimits_{t=0}^{T_{\mathrm{RS}}}\sum\nolimits_{i=1}^{n}
  \bigl\{(-C^{\mathrm{H_{2}}}q_{t}^{\mathrm{H_{2}}})\Delta t_{\mathrm{RS}}+C_{\mathrm{up,h}}^{\mathrm{EL}}\mu_{i,t}^{\mathrm{up,h}}+C_{\mathrm{up,c}}^{\mathrm{EL}}\mu_{i,t}^{\mathrm{up,c}}+C_{\mathrm{dow}n}^{\mathrm{EL}}\mu_{i,t}^{\mathrm{down}}+
  s(\sum\nolimits_{l=1}^{m}{P}_{l,t}^{\mathrm{WT,cut}}+{P}_{t}^{\mathrm{PV,cut}})\Delta t_{\mathrm{RS}}\bigr\},
\end{equation}
where $\boldsymbol{x}$ represents the decision variables, including the operational states and load power of ELs;
$m$ and $n$ are the number of WTs and ELs in each independent OReP2H sub-system, respectively;
$C^{\mathrm{H_{2}}}q_{t}^{\mathrm{H_{2}}}\Delta t_{\mathrm{RS}}$ is the revenue of selling hydrogen;
$C_{\mathrm{up,h}}^{\mathrm{EL}}\mu_{i,t}^{\mathrm{up,h}}+C_{\mathrm{up,c}}^{\mathrm{EL}}\mu_{i,t}^{\mathrm{up,c}}+C_{\mathrm{dow}n}^{\mathrm{EL}}\mu_{i,t}^{\mathrm{down}}$ represent the cost of state transition of an EL;
$s(\sum_{l=1}^{m}{P}_{l,t}^{\mathrm{WT,cut}}+{P}_{t}^{\mathrm{PV,cut}})\Delta t_{\mathrm{RS}}$ represents the penalty for curtailing renewable power. {\color{black} Detailed definitions of all variables and parameters are provided in the Nomenclature section (see~\ref{sec:appB}).}
%$n$ is the number of ELs;
%$m$ is the number of WTs;
%$\hat{P}_{l,t}^{\mathrm{WT,cut}}$ and $\hat{P}_{t}^{\mathrm{PV,cut}}$ indicate curtailed power of WTs and the PV plant, respectively;
%$\Delta t_{\mathrm{RS}}$ is the time step of scheduling;
%$\mu_{i,t}^{\mathrm{up,h}}$, $\mu_{i,t}^{\mathrm{up,c}}$, and $\mu_{i,t}^{\mathrm{down}}$ are binary variables indicating hot start-up, cold start-up, and shut-down of an EL.

The state transition and constraints of the electrolytic load of the ELs under different states can be expressed as follows \cite{qiu2023extended,zeng2024coordinated,yu2023optimal}:
\begin{gather}
  \mu_{i,t}^\mathrm{st}+\mu_{i,t}^\mathrm{sb}+\mu_{i,t}^\mathrm{sp} = 1,  \label{eq-2} \\
  \mu_{i,t}^\mathrm{down} = \mu_{i,t-1}^\mathrm{st}\mu_{i,t}^\mathrm{sp}, \label{eq-3} \\
  \mu_{i,t}^\mathrm{up,h}  =\mu_{i,t-1}^\mathrm{sb}\mu_{i,t}^\mathrm{st}, \label{eq-4} \\
  \mu_{i,t}^\mathrm{up,c}  =\mu_{i,t-1}^\mathrm{sp}\mu_{i,t}^\mathrm{st}, \label{eq-5} \\
  T_{\mathrm{min}}^\text{down}{ \mu _ { i , t }^{\mathrm{down}}}  \leq \sum\nolimits_{t}^{t-1+T_{\text{min}}^{\mathrm{down}}} {\mu_{i,t}^{\mathrm{sp}}}, \ \forall t\leq T_{\mathrm{RS}}+1-T_{\mathrm{min}}^{\mathrm{down}}, \label{eq-6} \\
  \mu_{i,t}^\mathrm{sb}P_\mathrm{sb}^\mathrm{EL}+\mu_{i,t}^\mathrm{st}r_{\mathrm{min}}^\mathrm{EL}S_i^\mathrm{EL} \leq P_{i,t}^\mathrm{EL,RS}\leq\mu_{i,t}^\mathrm{st}r_{\mathrm{max}}^\mathrm{EL}S_i^\mathrm{EL}+\mu_{i,t}^\mathrm{sb}P_\mathrm{sb}^\mathrm{EL}. \label{eq-7}
\end{gather}
%\begin{equation}\label{eq-7}
%(1-\mu_{i,t}^\mathrm{sb})P_{i,t}^\mathrm{AEL,RS}=k_i^\mathrm{H_2}q_{i,t}^\mathrm{H_2},
%\end{equation}
\noindent
where
%(\ref{eq-2})--(\ref{eq-5}) describe the states switching of AELs;
(\ref{eq-2})--(\ref{eq-5}) describe the transitions among  production, standby, and shut-down states of the ELs;
%$\mu_{i,t}^\mathrm{st}$, $\mu_{i,t1}^\mathrm{sb}$, and $\mu_{i,t}^\mathrm{sp}$ are binary variables indicating the production, standby, and shut-down states of the ELs, respectively;
(\ref{eq-6}) represents the minimal time interval between shut-down and start-up, where
$T_{\mathrm{min}}^\text{down}$ is set at 1 hour;
(\ref{eq-7}) constrains the load range and ramping limits of ELs under production and standby states.
%$P_{i,t}^\mathrm{EL,RS}$ is the load of the ELs;
%$P_\mathrm{sb}^\mathrm{EL}$ is the power consumption in standby state; $S_i^\mathrm{EL}$ is the capacity of a single EL;
%$r_{\mathrm{min}}^\mathrm{EL}$ and $r_{\mathrm{max}}^\mathrm{EL}$ are ramping limits. %, typically with $r_{\mathrm{min}}^\mathrm{EL}=10\%$ and $r_{\mathrm{max}}^\mathrm{EL} =120\%$.

The relationship between hydrogen production and electricity consumption is nonlinear and is influenced by the temperature and pressure of ELs. In this section, to ensure real-time performance of the RS, the relationship between the electrolytic power $P_{i,t}^\mathrm{EL,RS}$ and hydrogen production $q_{i,t}^\mathrm{H_2}$ is described using a fixed electricity consumption coefficient $k_i^\mathrm{H_2}$ following the typical approach used in the {\color{black}previous} work \cite{marocco2020study,marocco2022role,marocco2021optimal,oyewole2024optimal,viole2023renewable}:
\begin{equation}\label{eq-8}
  (1-\mu_{i,t}^\mathrm{sb})P_{i,t}^\mathrm{EL,RS}=k_i^\mathrm{H_2}q_{i,t}^\mathrm{H_2}.
\end{equation}

It is important to  note that different from the use of a linear model described here, a nonlinear curve of hydrogen production and electricity consumption derived from a commercial AEL produced by \emph{Peric Hydrogen Technologies} is employed to ensure accuracy for BESS size optimization and LCOH evaluation in Section \ref{sec:BatOptim}.
	
In addition to the ELs, energy shifting regulation is also provided by the BESS. The constraints for the BESS are listed as follows:
\begin{gather}
  SOC_{t}^{\mathrm{BES}}=SOC_{0}^{\mathrm{BES}}+\sum\nolimits_{i=t}^{T_{\mathrm{RS}}}(\eta^{\mathrm{BES}}P_{t}^{\mathrm{BES,C}}-{P_{t}^{\mathrm{BES,D}}}/{\eta^{\mathrm{BES}}})\Delta t_{\mathrm{RS}}, \label{eq-9}\\
  \left|SOC_{T_{\mathrm{RS}}}^{\mathrm{BES}}-SOC_{0}^{\mathrm{BES}}\right|\leq 0.05,   \label{eq-10} \\
  SOC_{\mathrm{min}}^{\mathrm{BES}}\leq SOC_{t}^{\mathrm{BES}} \leq SOC_{\mathrm{max}}^{\mathrm{BES}} \label{eq-11} \\
  0\leq P_t^{\mathrm{BES,C}}\leq\mu_t^{\mathrm{BES,C}}P_t^{\mathrm{BES,Cmax}}, \label{eq-12} \\
  0\leq P_t^{\mathrm{BES,D}}\leq\mu_t^{\mathrm{BES,D}}P_t^{\mathrm{BES,Dmax}}, \label{eq-13} \\
  0\leq\mu_t^{\mathrm{BES,C}}+\mu_t^{\mathrm{BES,D}}\leq1, \label{eq-14}
\end{gather}
\noindent
where (\ref{eq-9}) is the state of charge (SOC) of the BESS;
(\ref{eq-10})--(\ref{eq-11}) constraint the SOC variation range;
%$\eta^{\mathrm{BES}}$ is charging/discharging efficiency; $P_{t}^{\mathrm{BES,C}}$ and $P_{t}^{\mathrm{BES,D}}$ represent the charging and discharging power, respectively;
the final constraints (\ref{eq-12})--(\ref{eq-14}) prevent the BESS from charging and discharging simultaneously;
%and $\mu_t^{\mathrm{BES,C}}$ and $\mu_t^{\mathrm{BES,D}}$ are binary variables.

The balance of power for each time step in the RS is finally required, and is described by:
\begin{gather}
  {\color{black}\sum\nolimits_{l=1}^m\widehat{P}_{l,t}^{\mathrm{WT}} + \widehat{P}_t^{\mathrm{PV}} - \sum\nolimits_{l=1}^m{P}_{l,t}^{\mathrm{WT,cut}} -
  {P}_t^{\mathrm{PV,cut}}+P_t^{\mathrm{BES,D}}=P_t^{\mathrm{BES,C}}+\sum\nolimits_{i=1}^nP_{i,t}^{\mathrm{EL,RS}}}, \label{eq-15}\\
  0\leq{P}_{l,t}^\text{WT,cut}\leq\widehat{P}_{l,t}^\text{WT}, \ \text{for} \ l = 1,\ldots,m, \label{eq-16} \\
  0 \leq{P}_t^{\text{PV,cut}} \leq \widehat { P }_t^{\text{PV}}, \label{eq-17}
\end{gather}
\noindent
where (\ref{eq-16})--(\ref{eq-17}) represent the allowable power curtailments of WTs and PV under RS stage.

%{\color{black}${P}_{l,t}^{\mathrm{WT,cut}}$} and {\color{black}${P}_t^{\mathrm{PV,cut}}$} are the power curtailments of WT and PV, respectively.

Finally, we summarize (\ref{eq-1})--(\ref{eq-17}) to establish the MILP-based RS model. This RS model is similar to those  used in the previous studies for ReP2HS planning \cite{oyewole2024optimal,wang2023optimising,ibagon2023techno,yang2020planning,shao2023risk}, and it can be efficiently solved in polynomial time since the computational complexity is $O(2^{1248} \cdot $poly$(4, 48)$, where $ $poly$(4, 48)$ denotes the polynomial-time complexity of solving the linear programming subproblem with fixed integer variables. {\color{black}Table~\ref{tab:AvgCompuTime} further verifies the computational efficiency, showing that even with 16 ELs in a single subsystem of OReP2HS, the rolling scheduling can be solved by \emph{Gorubi 12.0} within 7.434~seconds on average, fully meeting the computational requirement.} However, while prior works only consider intra-day scheduling (as discussed in Section \ref{sec:review}), we further investigate the influence of seconds-level load control of the AELs (Section \ref{sec:SFL}) and the grid-forming control of the BESS at the millisecond scale (Section \ref{sec:GFM}) on the sizing of the BESS.
\begin{table}[h]\footnotesize
  \renewcommand{\arraystretch}{1.15}
  \caption{{\color{black}Average computation time of rolling scheduling under different number of ELs.}}
  \label{tab:AvgCompuTime}
  \vspace{0pt}
  \centering
  \begin{tabular}{cc}
  \hline \hline
   {\color{black} Number of Electrolyzers} & {\color{black} Average Computation Time of 10 Scheduling Runs (s)} \\
  \hline
   {\color{black} 4} & {\color{black} 0.375}  \\
   {\color{black} 8} & {\color{black} 0.964}  \\
   {\color{black} 12} & {\color{black} 2.502} \\
   {\color{black} 16} & {\color{black} 7.434} \\
  \hline \hline
  \end{tabular}
\end{table}

\subsection{Load Following}
\label{sec:SFL}
%The LF module adjusts the load of each EL in real time. It aims to optimize the load regulation capacity of ELs to balance power within seconds. This module includes \emph{fast forecasting}, \emph{load following control}, \emph{SOC correction}, and \emph{load allocation}, as shown in Fig. \ref{fig:EMS}.
%
%The fast forecasting sub-module predicts renewable power for the next few seconds (determined by $\Delta t_\mathrm{LF}$) using a discrete moving average (MA) process, as shown in (\ref{eq-18}). It uses real-time active power measurement from WTs and the PV plant to generate the load correction command $P_t^\mathrm{LF}$ according to (\ref{eq-19}).
%\begin{gather}
%  \hat{P}_t^\mathrm{RES}=\alpha\hat{P}_{t-\Delta t_\mathrm{LF}}^\mathrm{RES}+(1-\alpha)\frac1{q_\mathrm{LF}}\sum\nolimits_{k=1}^{q_\mathrm{LF}}\bigl(\sum\nolimits_{l=1}^mP_{l,t-k\Delta t_\mathrm{LF}}^\mathrm{WT}+P_{t-k\Delta t_\mathrm{LF}}^\mathrm{PV}\bigr), \label{eq-18} \\
%  P_t^\mathrm{LF}=\hat{P}_t^\mathrm{RES}-\sum\nolimits_{i=1}^nP_{i,t-\Delta t_\mathrm{LF}}^\mathrm{EL}, \label{eq-19}
%\end{gather}
%\noindent
%where $\hat{P}_t^\mathrm{RES}$ is the forecasted renewable power; $q_\mathrm{LF}$ and $\alpha$ are the order and smoothing coefficient of the MA process, set at 4 and 0.6, respectively; $\Delta t_\mathrm{LF}$ is the LF time step, typically on the order of  seconds (the effect of $\Delta t_\mathrm{LF}$ on the EMS performance is discussed in Section \ref{sec:Sensitivity}).

The LF module adjusts the load of each EL in real time, {\color{black}serving as a rolling correction mechanism to handle renewable power uncertainties. It aims to optimize the load regulation capability of the ELs to balance renewable power fluctuations within seconds.} This module includes \emph{fast forecasting}, \emph{load following control}, \emph{SOC correction}, and \emph{load allocation}, as shown in Fig. \ref{fig:EMS}.

To clearly distinguish time resolutions in our model, we denote \( \tau \) as the time index for the LF module (typically every 5 seconds), while \( t \) continues to denote the 5-min time step used in the RS module.

Since the LF control operates on a seconds time scale, faster power forecasting is required to meet its high-frequency control demands {\color{black}and to capture short-term fluctuations of renewable generation}. Therefore, a rapid power forecasting module based on a discrete moving-average (MA) process is integrated into the LF strategy. It uses real-time active power measurement from WTs and the PV plant to generate the load correction command \( P_{t,\tau}^\mathrm{LF} \), as shown in (\ref{eq-18})--(\ref{eq-19}).
\begin{gather}
  \tilde{P}_\tau^\mathrm{RES}=\alpha\tilde{P}_{\tau-1}^\mathrm{RES}+(1-\alpha)\frac1{q_\mathrm{LF}}\sum\nolimits_{k=1}^{q_\mathrm{LF}}\bigl(\sum\nolimits_{l=1}^mP_{l,\tau-k}^\mathrm{WT}
  +P_{\tau-k}^\mathrm{PV}\bigr), \label{eq-18} \\
  P_{t,\tau}^\mathrm{LF}=\tilde{P}_\tau^\mathrm{RES}-\sum\nolimits_{i=1}^nP_{i,t}^\mathrm{EL,RS}, \label{eq-19}
\end{gather}

\noindent
where \( \tilde{P}_\tau^\mathrm{RES} \) is the forecasted renewable power at time index $\tau$ used in LF; \( q_\mathrm{LF} \) and \( \alpha \) are the order and smoothing coefficient of the MA process, set at 4 and 0.6, respectively; subscription \( (t,\tau) \) denote the LF command at seconds-level index $\tau$ within the 5-minute-level RS period $t$. The effect of different seconds-level time steps of LF on EMS performance is discussed in Section~\ref{sec:Sensitivity}.

The load-following sub-module tracks $ P_{t,\tau}^\mathrm{LF} $ using a PI controller, followed by the use of a SOC correction control sub-module to prevent battery over-charging/discharging. The total correction $P_{t,\tau}^\mathrm{LF,ref}$ is given by:
\begin{gather}\label{eq-20}
  P_{t,\tau}^\mathrm{LF,SOC}=(K_{\mathrm{SOC}}P_{t,\tau}^\mathrm{LF}+\beta), \\
  P_{t,\tau}^{\mathrm{LF,ref}}=K_P(P_{t,\tau}^\mathrm{LF,SOC}+\sum\nolimits_{i=1}^nP_{i,t}^\mathrm{EL,RS}-\sum\nolimits_{i=1}^nP_{i,t,\tau}^\mathrm{EL})+\int_0^t K_I(P_{t,\tau}^\mathrm{LF,SOC}+\sum\nolimits_{i=1}^nP_{i,t}^\mathrm{EL,RS}-\sum\nolimits_{i=1}^nP_{i,t,\tau}^\mathrm{EL}) \mathrm{d}t,
\end{gather}
where $K_{\mathrm{SOC}}$ and $\beta$ are constant coefficients, and $K_P$ and $K_I$ are the proportional and integral gains of the PI controller, respectively. $P_{i,t,\tau}^\mathrm{EL}$ is the actual load power of $i$-th EL at LF index $\tau$ within RS period $t$.

The allocation of $P_{t,\tau}^{\mathrm{LF,ref}}$ to each EL is based on their adjustable capacity, as shown in (\ref{eq-21}). If $P_{t,\tau}^{\mathrm{LF,ref}}\geq0$, indicating a need to increase the electrolytic load, it is distributed based on the upward regulation capacities of ELs, and vice versa for a decrease in the load.
\begin{equation}\label{eq-21}
  P_{i,t,\tau}^\text{EL,LF}=
  \begin{cases}
  	 \dfrac{S_{i}^{\mathrm{EL}}-P_{i,t,\tau}^{\mathrm{EL}}}{\sum_{i=1}^{n}(S_{i}^{\mathrm{EL}}-P_{i,t,\tau}^{\mathrm{EL}})} \times P_{t,\tau}^{\mathrm{LF,ref}},\
  \text{for} \ P_t^\mathrm{LF,ref}\geq0, \\
  \dfrac{-P_{i,,\tau}^{\mathrm{EL}}}{\sum_{i=1}^{n}P_{i,t,\tau}^{\mathrm{EL}}}  \times P_{t,\tau}^{\mathrm{LF,ref}},\ \text{for} \  P_{t,\tau}^\mathrm{LF,ref}<0.
  \end{cases}
\end{equation}

Finally, the load reference sent to each EL including LF and RS commands are described by:
\begin{equation}\label{eq-22}
P_{i,t,\tau}^\mathrm{EL}=P_{i,t}^\mathrm{EL,RS}+P_{i,,\tau}^\mathrm{EL,LF},
\end{equation}
where $P_{i,t}^\mathrm{EL,RS}$ is the baseline load of the $i$-th EL determined by the RS module presented in Section \ref{sec:RS} that is updated every $\Delta t_\mathrm{RS}$; and $P_{t,\tau}^{\mathrm{LF,ref}}$ is the real-time correction command given by (\ref{eq-21}) that is updated typically every 5 seconds {\color{black}which effectively mitigate the uncertainty of renewable power outputs.}

\subsection{Emergency Handling}
\label{sec:SCODE}

In emergencies, such as the tripping of a WT, PV plant, or EL, a common \textit{N}$-1$ fault in power systems, the stability of the OReP2HS relies heavily on the grid-forming BESS for rapid power support and prevention of unnecessary blackouts. However, continuous reliance on the BESS during these emergencies can deplete its energy reserves, increasing the required BESS capacity and, consequently, the investment costs. %Under the V/F droop-based grid-forming control, active power and frequency are coupled, allowing for the estimation of imbalanced power through system frequency.
To address this, we propose an EH strategy that coordinates WTs, PV, and ELs to reduce the energy demand on the BESS during emergencies.

The procedure for EH is outlined below. First, it estimates the imbalanced power through  system frequency and then regulates the BESS along with PV/WTs/ELs. The severity of the emergency is assessed using the rate of change of frequency (RoCoF) and maximal frequency deviation. When EH is activated, a power shedding command is dispatched to the renewable power sources or electrolytic loads, depending on whether the power loss originated from the load or source side.The ELs or generators respond to the power shedding command with maximal ramping, reducing the power support demand from the BESS. This facilitates rapid restoration of power balance within the OReP2HS, potentially reducing BESS capacity.

{\color{black}To ensure practical applicability, the proposed EH strategy adopts a lookup-table-based activation approach (see Table~3, tailored for the system shown in Fig.~2). Although it may not be as advanced as adaptive methods, it remains one of the most practical and straightforward options for real-world implementation of OReP2HS and the simulation results presented in Section~5.1.2 demonstrate its effectiveness.}
%,  with activation criteria set by simulation tests.

Notably, for EH activation, both RoCoF and maximal frequency deviation must fall within predefined intervals specified in Table \ref{tab:Emergency}. If either indicator falls outside these intervals, EH cannot be activated; if both fall outside all intervals, EH is deactivated.

\begin{table}[tb]\footnotesize
  \renewcommand{\arraystretch}{1.8}
  \caption{Load/Generation shedding command and activation criteria of the EH}
  \label{tab:Emergency}
   \vspace{-6pt}
  %\newcolumntype{C}{>{\centering\arraybackslash}X}
  \centering
  \begin{tabular}{ccc}
  \hline \hline
  \vspace{-30pt}
  \\ \vspace{-1.5pt}
   {\tabincell{c}{Interval}\vspace{1.5pt}}   &
   {\tabincell{c}{Load/Generation \\shedding}}   &
   {\tabincell{c}{Activation criteria (Hz and Hz/s)}}    \\
   \hline
   {\tabincell{c}{SE6}}   &  {\tabincell{c}{6 MW}}   &  {\tabincell{c}{$f_t^{\mathrm{PCC}}\leq49.65 \enspace \mathrm{and} \enspace \mathrm{d}f_t^{\mathrm{PCC}}/\mathrm{d}t\leq-11.5$\\
   $f_t^{\mathrm{PCC}}\geq50.35 \enspace \mathrm{and} \enspace \mathrm{d}f_t^{\mathrm{PCC}}/\mathrm{d}t\geq11.5 $}}   \\

  {\tabincell{c}{SE5}}   &  {\tabincell{c}{5 MW}}   &  {\tabincell{c}{$49.70\leq f_t^{\mathrm{PCC}}<49.65 \enspace \mathrm{and} \enspace -11.5<\mathrm{d}f_t^{\mathrm{PCC}}/\mathrm{d}t\leq-10.0$ \\ $50.30\leq f_t^{\mathrm{PCC}}<50.35 \enspace \mathrm{and} \enspace 10.0\leq\mathrm{d}f_t^{\mathrm{PCC}}/\mathrm{d}t<11.5$}}   \\

  {\tabincell{c}{SE4}}   &  {\tabincell{c}{4 MW}}   &  {\tabincell{c}{$49.75\leq f_t^{\mathrm{PCC}}<49.70 \enspace \mathrm{and} \enspace -10.0<\mathrm{d}f_t^{\mathrm{PCC}}/\mathrm{d}t\leq-8.0$ \\ $50.25\leq f_t^{\mathrm{PCC}}<50.30 \enspace \mathrm{and} \enspace 8.0\leq\mathrm{d}f_t^{\mathrm{PCC}}/\mathrm{d}t<10.0$}}   \\

  {\tabincell{c}{SE3}}   & {\tabincell{c}{3 MW}}   &  {\tabincell{c}{$49.80\leq f_t^{\mathrm{PCC}}<49.75 \enspace \mathrm{and} \enspace -8.0<\mathrm{d}f_t^{\mathrm{PCC}}/\mathrm{d}t\leq-6.5$ \\ $50.20\leq f_t^{\mathrm{PCC}}<50.25 \enspace \mathrm{and} \enspace 6.5\leq\mathrm{d}f_t^{\mathrm{PCC}}/\mathrm{d}t<8.0$}}   \\

  {\tabincell{c}{SE2}}   &  {\tabincell{c}{2 MW}}   &  {\tabincell{c}{$49.85\leq f_t^{\mathrm{PCC}}<49.80 \enspace \mathrm{and} \enspace -6.5<\mathrm{d}f_t^{\mathrm{PCC}}/\mathrm{d}t\leq-4.5$ \\ $50.15\leq f_t^{\mathrm{PCC}}<50.20 \enspace \mathrm{and} \enspace 4.5\leq\mathrm{d}f_t^{\mathrm{PCC}}/\mathrm{d}t<6.5$}}   \\

  {\tabincell{c}{SE1}}   & {\tabincell{c}{1 MW}}   &  {\tabincell{c}{$49.85\leq f_t^{\mathrm{PCC}}<49.80 \enspace \mathrm{and} \enspace -4.5<\mathrm{d}f_t^{\mathrm{PCC}}/\mathrm{d}t\leq-2.5$ \\ $50.15\leq f_t^{\mathrm{PCC}}<50.20 \enspace \mathrm{and} \enspace 2.5\leq\mathrm{d}f_t^{\mathrm{PCC}}/\mathrm{d}t<4.5$}}\vspace{3pt} \\

  \hline \hline
  \end{tabular}
\end{table}

\subsection{V/F Droop-based Grid-forming Control}
\label{sec:GFM}

The V/F droop-based grid-forming control is integrated into the built-in controller of the BESS. It enables the BESS to actively provide frequency and voltage references to the OReP2HS, and its operating principle is described in (\ref{eq-23})--(\ref{eq-24}) \cite{rosso2021gridforming, lasseter2019gridforming}.
\begin{gather}
  \omega_{ t_{\mathrm{sim}}}=\omega_{ t_{\mathrm{sim}}}^{\mathrm{ref}}+K_f(P_{ t_{\mathrm{sim}}}^{\mathrm{BES,ref}}-P_{ t_{\mathrm{sim}}}^{\mathrm{BES}}), \label{eq-23}\\
  V_{ t_{\mathrm{sim}}}=V_{ t_{\mathrm{sim}}}^{\mathrm{ref}}+K_v(Q_{ t_{\mathrm{sim}}}^{\mathrm{BES,ref}}-Q_{ t_{\mathrm{sim}}}^{\mathrm{BES}}), \label{eq-24}
\end{gather}
where ${t_{\text{sim}}}$ is the time of the electromagnetic transient simulation, typically on the order of milliseconds; $\omega_t$ and $\omega_t^{\mathrm{ref}}$ are the angular frequency ($\omega_t = 2 \pi f_t$) and its reference; $V_t$ and $V_t^{\mathrm{ref}}$ are the voltage at the point of common connection (PCC) and its reference; and $K_f$ and $K_v$ are droop coefficients, representing the frequency and voltage variation w.r.t. the difference between the instantaneous outputs of the BESS  $P_t^{\mathrm{BES}}/Q_t^{\mathrm{BES}}$ and its setpoint  $P_t^{\mathrm{Bat,ref}}/Q_t^{\mathrm{Bat,ref}}$. A comprehensive introduction to grid-forming control can be found in  \cite{lasseter2019gridforming}.

However, implementation of grid-forming control in a BESS is a complex task that requires coordination among the outer V/F droop control, inner current control, measurement, filter, and PWM generator \cite{rosso2021gridforming}. This procedure involves interactions between different control loops on the millisecond time scale, necessitating the use of an electromagnetic transient simulation model, as described in detail in Sections \ref{sec:Model} and \ref{sec:appA}.

\begin{figure}[H]
\centering
\includegraphics[width=5.8in]{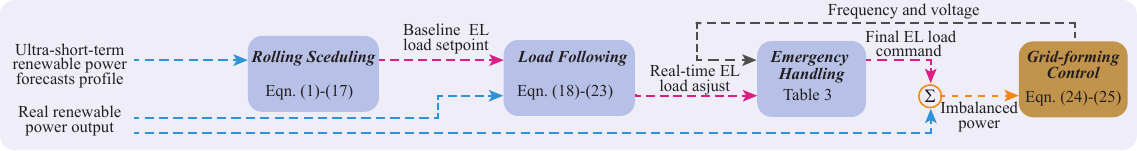}
%\vspace{-6pt}
\caption{\color{black} Schematic diagram of the inputs/outputs interactions among the different modules of the proposed multi-timescale EMS.}
\label{fig:SimpleEMS}
\end{figure}

{\color{black}For clarity, the inputs/outputs interactions among the different modules of the proposed multi-timescale EMS are illustrated in Fig.~\ref{fig:SimpleEMS}. The \textit{Rolling Scheduling} module determines the baseline EL load setpoints every 4 hours based on ultra-short-term renewable forecasts and system operational constraints~(\ref{eq-2})--(\ref{eq-17}). The \textit{Load Following} module refines these setpoints in real time using actual renewable power outputs to generate fast load adjustment commands. The \textit{Emergency Handling} module monitors frequency and voltage deviations; once abnormal conditions are detected, it issues corrective control signals follow  Table~\ref{tab:Emergency} to prevent system instability. The resulting power imbalance of renewable generations and ELs is then compensated by the \textit{Grid-forming Control} of the BESS, which regulates system frequency and voltage through transient power support. This structure enables coordinated operation across timescales ranging from milliseconds to hours, ensuring both transient stability and long-term energy balance of the OReP2HS.}

\section{Optimal Size of the BESS}
\label{sec:BatOptim}

\begin{landscape}
\begin{figure}[tb]
  \centering
  \includegraphics[width=9.2 in]{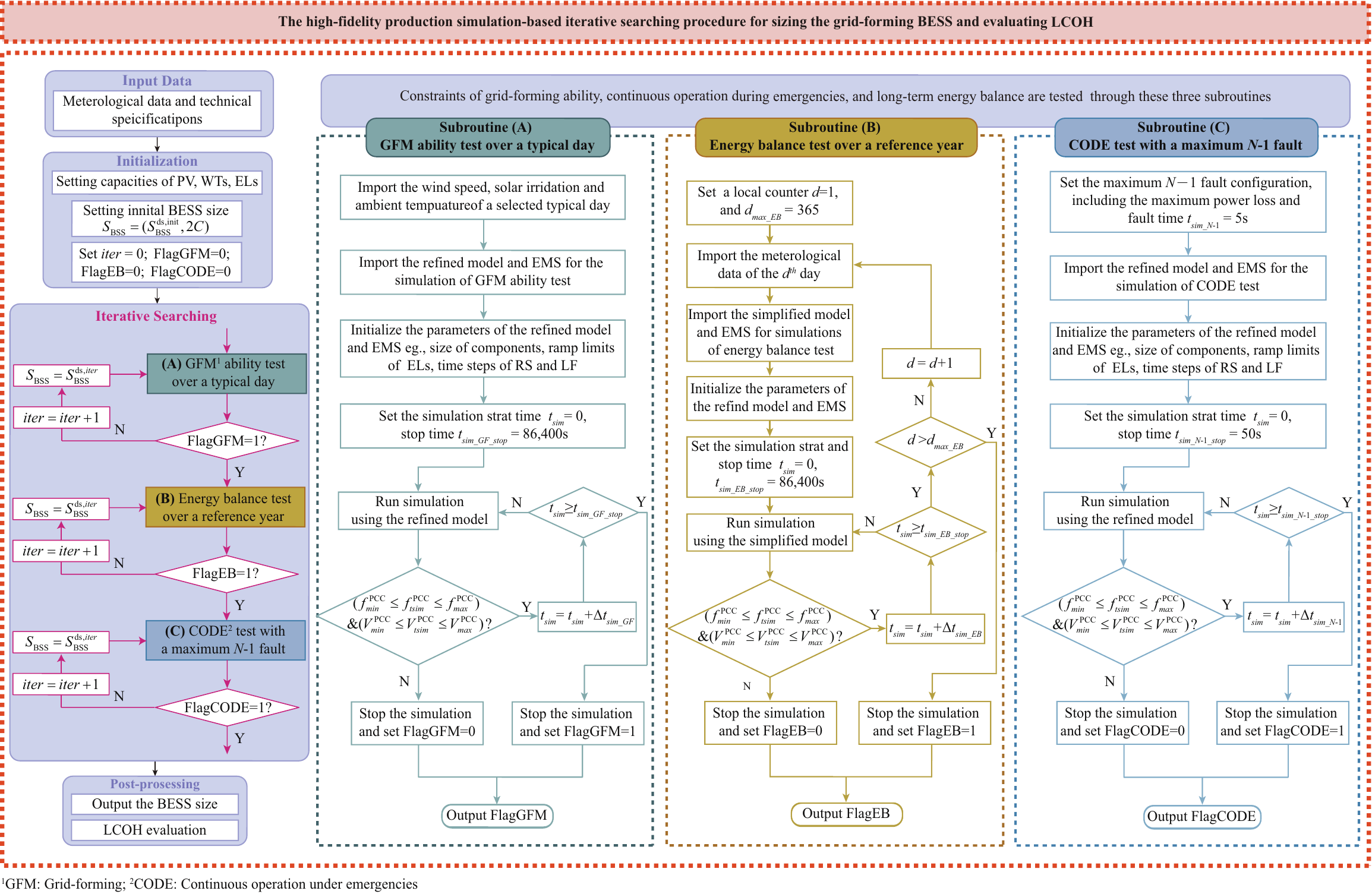}
  \caption{Framework of the iterative search procedure for sizing of the grid-forming BESS.}
  \label{fig:Solving}
\end{figure}
\end{landscape}

As mentioned in Section \ref{sec:Requorement}, sizing the grid-forming BESS involves testing the grid-forming capability, continuous operation during emergencies, and long-term energy balance. These constraints involve complex interactions between WTs, PV plant, ELs, and the BESS, making explicit quantification extremely difficult. Existing methods, such as those in \cite{oyewole2024optimal, viole2023renewable, wang2023optimising, ibagon2023techno, yang2020planning, shao2023risk} that use convex optimization to find the optimal BESS size, are not applicable.
Therefore, we employ a high-fidelity production simulation-based iterative searching procedure, based on the proposed multi-timescale EMS, to evaluate the most cost-effective BESS size. The iterative search procedure is illustrated in Fig. \ref{fig:Solving}, with the objective and constraints described below.

\subsection{Searching Objective}
\label{sec:BatObj}

The objective of the iterative search is to minimize the LCOH, as described by:
\begin{equation}\label{eq-25}
  \min\ LCOH=\frac{(C_\mathrm{init} + C_{\mathrm{O\&M}}^\mathrm{fixed} + C_{\mathrm{O\&M}}^\mathrm{vari} ) }{M^\mathrm{H_{2}}},
\end{equation}
\noindent
where $C_\mathrm{init}$ is the annualized capital cost; $C_{\mathrm{O\&M}}^\mathrm{fixed}$ is the fixed annual operation and maintenance (O$\&$M) cost; $C_{\mathrm{O\&M}}^\mathrm{vari}$ is the variable annual O$\&$M cost; and $M^\mathrm{H_{2}}$ is the annual hydrogen yield.

The annualized capital cost $C_\mathrm{init}$ is calculated as follows:
\begin{gather}
  C_{\mathrm{init}}=\sum\nolimits_{j\in\Omega_{\mathrm{F}}}CRF(r,T_j^{\mathrm{LCC}})S_jc_j^{\mathrm{unit}},\ \Omega_{\mathrm{F}}=\{\mathrm{WT,PV,EL,BES,TL}\},
  \label{eq-26} \\
  CRF(r,T_{j}^{\mathrm{LCC}})=\frac{r(1+r)^{T_{j}^{\mathrm{LCC}}}}{(1+r)^{T_{j}^{\mathrm{LCC}}-1}}, \ \text{for}\ j\in\Omega_{\mathrm{F}}, \label{eq-27}
\end{gather}
where $\mathrm{TL}$ represents transmission lines; $S_j$ is the device capacity; $c_j^{\mathrm{unit}}$ is the unit cost of device; $CRF(r,T_{j}^{\mathrm{LCC}})$ is the capital recovery factor; $r$ is the discount rate (set at 8\%); and $T_{j}^{\mathrm{LCC}}$ is the lifetime of the facilities.

The fixed annual O\&M cost $C_{\mathrm{O\&M}}^\mathrm{fixed}$ is assumed to be a fixed proportion of the initial investment, as expressed by:
\begin{equation}\label{eq-28}
  C_\mathrm{O\&M}^\mathrm{fixed}=\sum\nolimits_{j\in\Omega_{\mathrm{F}}}\lambda_j^\mathrm{O\&M}S_jc_j^\mathrm{unit},
\end{equation}
where $\lambda_j^\mathrm{O\&M}$ represent the fixed proportion.

The variable annual O$\&$M cost $C_{\mathrm{O\&M}}^\mathrm{vari}$ comprises the BESS replacement cost $C_{\mathrm{BES}}^{\text{rep}}$, and the benefit of recycling the retired battery $C_{\mathrm{BES}}^{\text{rec}}$, written as:
\begin{gather}
  C_{\mathrm{O\&M}}^{\mathrm{vari}}=C_{\mathrm{BES}}^{\mathrm{rep}}-C_{\mathrm{BES}}^{\mathrm{rec}}, \label{eq-29}\\
  C_{\mathrm{BES}}^{\text{rep}} = \sum\nolimits _ { k = 1 }^{K_{\mathrm{BES}}^{\text{rep}} }S_{\mathrm{BES}}^{\text{rep}} c _ {\mathrm{BES}} ^ {\text{rep}} ( 1 + r ) ^ { - k T _ {\mathrm{BES}} ^ {\text{LCC}} / ( K _ {\mathrm{BES}} ^\text{rep}{ + 1 })}, \label{eq-30} \\
  C_{\mathrm{BES}}^{\text{rec}} = \sum\nolimits _ { k = 1 }^{K_{\mathrm{BES}}^{\text{rep}} }S_{\mathrm{BES}}^{\text{rep}} c _ {\mathrm{BES}} ^ {\text{rec}} ( 1 + r ) ^ { - k T _ {\mathrm{BES}} ^ {\text{LCC}} / ( K _ {\mathrm{BES}}^\text{rep}{ + 1 })}, \label{eq-31}
\end{gather}
\noindent
where $S_{\mathrm{BES}}^{\text{rep}}$ is the capacity of the retired battery; $c _ {\mathrm{BES}} ^ {\text{rep}}$ and $c _ {\mathrm{BES}} ^ {\text{rec}}$ are the unit cost and revenue of battery replacement and recycling, respectively; and $K _ {\mathrm{BES}} ^\text{rep}$ is the number of battery replacements.

The number of battery replacements $K _ {\mathrm{BES}} ^\text{rep}$ due to degradation is estimated according to:
\begin{gather}
K_{\mathrm{BES}}^{\mathrm{rep}}= \left[ {T_{\mathrm{BES}}^{\mathrm{LCC}}}/{T_{\mathrm{BES}}^{\mathrm{Fail}}}-1 \right], \label{eq-32} \\
T_{\mathrm{BES}}^{\text{Fail}} = { \Delta E _ {\mathrm{BES}} ^ {\text{degra,}\text{max}}}/{ \Delta E _ {\mathrm{BES}} ^ {\text{degra,year}}}, \label{eq-33}
\end{gather}
where $[ \cdot ]$ is the ceiling function; $ \Delta E _ {\mathrm{BES}} ^ {\text{degra,year}}$ is the annual battery degradation, which is estimated based on the methods in \cite{debnath2014quantifying,lam2012practical} from the simulation data obtained by the high-fidelity model presented in Section \ref{sec:Model}; and $\Delta E _ {\mathrm{BES}} ^ {\text{degra,}max}$ is the maximum allowable degradation and is typically 20\%.

The annual hydrogen yield is calculated as:
\begin{equation}\label{eq-34}
  M^{\mathrm{H}_2}=\sum\nolimits_{t=0}^{8760/\Delta t_{\text{sim}}} \sum\nolimits_{i=1}^{n}(\alpha_{\mathrm{P2H}}(q_t^{\mathrm{H}_2},P_{t}^{\mathrm{EL}})P_{i,t}^{\mathrm{EL}}\Delta t_{\text{sim}},
\end{equation}
where $\alpha_{\mathrm{P2H}}(q_t^{\mathrm{H}_2},P_{t}^{\mathrm{EL}})$ is the conversion efficiency function obtained using a 5 MW-rated AEL produced by \emph{Peric Hydrogen Technologies}.

\subsection{Feasibility Testing Constraints}
\label{sec:BatConstrain}

The requirements listed in Section \ref{sec:Requorement} serve as constraints in the search for the optimal BESS size. They include grid-forming ability, continuous operation during emergencies, and long-term energy balance.
The primary goal of setting these constraints is to ensure frequency and voltage stability and maintain energy balance across timescales.
However, these constraints involve multi-timescale transient processes (detailedly discussed in Section \ref{sec:GFM}) and the complex interactions among WTs, PV, ELs, and BESS, making explicit quantification challenging.
To address this issue, we use refined and simplified electromagnetic transient simulation models presented in Section \ref{sec:Model} and conduct high-fidelity time-domain simulations to test whether these constraints are satisfied; the detailed procedure is given in Section \ref{sec:BatSolving}.

The frequency and voltage at PCC are key indicators of the OReP2HS stability\cite{mandal2020frequency,murray2021voltage}. Any imbalance in the power and energy in the OReP2HS will finally lead to frequency or voltage changes. Therefore, these indicators serve as criteria for testing all of the constraints:
\begin{equation}\label{eq-35}
  \begin{cases}
    f_{\text{min}}^\mathrm{PCC}\leq f_{t_{\text{sim}}}^\mathrm{PCC}\leq f_{\text{max}}^\mathrm{PCC},\\
    V_{\text{min}}^\mathrm{PCC}\leq V_{t_{\text{sim}}}^\mathrm{PCC}\leq V_{\text{max}}^\mathrm{PCC},
  \end{cases}
\end{equation}
where $f_{\text{min}}^\mathrm{PCC}$ and $f_{\text{max}}^\mathrm{PCC}$ are the minimum and maximum frequency limits, respectively; and $V_{\text{min}}^\mathrm{PCC}$ and $V_{\text{max}}^\mathrm{PCC}$ are the bus voltage limits; and $f_{\text{tsim}}^\mathrm{PCC}$ and $V_{\text{tsim}}^\mathrm{PCC}$ represent the frequency and voltage at each time step of the electromagnetic transient simulation.

\subsection{Procedure for Searching for the Optimal BESS Size}
\label{sec:BatSolving}

\begin{table}[tb]\footnotesize
  \renewcommand{\arraystretch}{1.25}
  \caption{BESS size preset table for iterative searching}
  \label{tab:Preset}
  \vspace{-6pt}
  %\newcolumntype{C}{>{\centering\arraybackslash}X}
  \centering
  \begin{tabular}{ccccccc}
  \hline \hline
  & \multicolumn{6}{c}{Preset of the BESS size (battery capacity, charging/discharging rate)} \\ %合并7列

  \hline

  \makebox[0.04\textwidth][c]{\tabincell{c}{Iterative \\variable}} &  $S_{\mathrm{BES}}^\mathrm{ds,0}$  & $S_{\mathrm{BES}}^\mathrm{ds,1}$  & $S_{\mathrm{BES}}^\mathrm{ds,2}$  & \makebox[0.01\textwidth][c]{$\cdots$}  & $S_{\mathrm{BES}}^{\mathrm{ds,}iter-1}$ & $S_{\mathrm{BES}}^{\mathrm{ds,}iter}$ \\

  \makebox[0.04\textwidth][c]{Size} & \makebox[0.07\textwidth][c]{\makecell[c]{$(S_{\mathrm{BES}}^\mathrm{init},2\mathrm{C})$}}
  & \makebox[0.15\textwidth][c]{$(S_{\mathrm{BES}}^\mathrm{init}+\Delta S^\mathrm{BES},2\mathrm{C})$}
  & \makebox[0.15\textwidth][c]{$(S_{\mathrm{BES}}^\mathrm{init}+\Delta S^\mathrm{BES},3\mathrm{C})$}  & \makebox[0.01\textwidth][c]{$\cdots$}
  &  \makebox[0.187\textwidth][c]{$(S_{\mathrm{BES}}^\mathrm{init}+[\frac{iter}{2}]\Delta S^\mathrm{BES},2\mathrm{C})$}
  &
  \makebox[0.187\textwidth][c]{$(S_{\mathrm{BES}}^\mathrm{init}+[\frac{iter}{2}]\Delta S^\mathrm{BES},3\mathrm{C})$} \\
  \hline \hline

  \multicolumn{7}{l}{$iter\in\mathbb{N}; \Delta S^\mathrm{BES}$ is the incremental step; $[ \cdot ]$ is the ceiling function. }   \\

  \end{tabular}
\end{table}

The optimal BESS size in the OReP2HS is evaluated by a high-fidelity-simulation-based iterative search, as illustrated in Fig. \ref{fig:Solving}.

The process begins by setting an initial BESS size, which includes two key parameters: battery capacity {\color{black}($S_{\mathrm{BES}}^\mathrm{init}$)} and charging/discharging rate {\color{black}(C).} Then, the iterative search begins. Three subroutines sequentially simulate and test the constraints of grid-forming ability, continuous operation during emergencies, and long-term energy balance. If any subroutine fails the test, the battery capacity is corrected according to the preset table (shown in Table \ref{tab:Preset}), and the subroutine is re-executed until it passes the test. {\color{black}In Table~\ref{tab:Preset}, the notation ($S_{\mathrm{BES}}^\mathrm{init}+\Delta S^\mathrm{BES}, 2\mathrm{C}$) represents a BESS configuration with a converter-to-capacity ratio of 2 C, meaning that the converter rated power is twice the battery energy capacity in MW per MWh.}

The features of each subroutine are designed to address the specific aspects of the BESS sizing process, using different simulation models presented in Section \ref{sec:Model}.
The grid-forming ability test simulates the daily operation of the OReP2HS using the proposed EMS and a refined model that accounts for the switching dynamics of the converters and the transient behaviors of PV, WTs, BESS, and ELs. The simulation employs a time step $\Delta t_{\text{sim}\_\text{GF}}$ of ten microseconds to capture these dynamics accurately.
The energy balance test uses the simplified model that considers the transient behaviors of PV, WTs, BESS, and ELs but omits the switching dynamics of the converters. It performs a production simulation over a reference year with a time step $\Delta t_{\text{sim}\_\text{EB}}$ of one hundred microseconds. The simulation is parallelized over multiple CPUs to speed up the computation.
The CODE test  evaluates the maximal \textit{N}$-1$ fault scenarios for both the source and load of the OReP2HS using the refined model.  The simulation time step $\Delta t_{\text{sim}\_N-1}$ for this test is the same as $\Delta t_{\text{sim}\_\text{GF}}$. Detailed flowcharts of each subroutine are presented in Fig. \ref{fig:Solving}.

Once all subroutines successfully meet the constraints, we record the feasible BESS sizes. The LCOH for each feasible size is then calculated using (\ref{eq-25})--(\ref{eq-34}), and the BESs size with the lowest LCOH is selected as optimal. Although the BESS sizing problem in (\ref{eq-25})--(\ref{eq-35}) is formulated as a single-objective model minimizing LCOH, it inherently considers multiple competing objectives. For instance, battery life are embedded in the O$\&$M costs (constraints (\ref{eq-29})--(\ref{eq-33})); the nonlinear efficiency of ELs, affecting economic efficiency, is captured in constraint (\ref{eq-34}); system stability in terms of voltage and frequency is enforced through constraint(\ref{eq-35}). Thus, the model can be also viewed as a multi-objective framework with a preference for economic efficiency.

\section{Results and Discussions}
\label{sec:CaseStudy}

The case study focuses on the OReP2HS planned in Inner Mongolia, China, as described in Section \ref{sec:Config}.
The simulation models, detailed in Section \ref{sec:Model}, are developed in \emph{Matlab/Simulink}. The size of the grid-forming BESS and LCOH are evaluated using the search procedure established in Section \ref{sec:BatOptim}.

\begin{figure}[tb]
	\centering
	\begin{subfigure}{0.45\linewidth}
		\centering
		\includegraphics[width=0.95\linewidth]{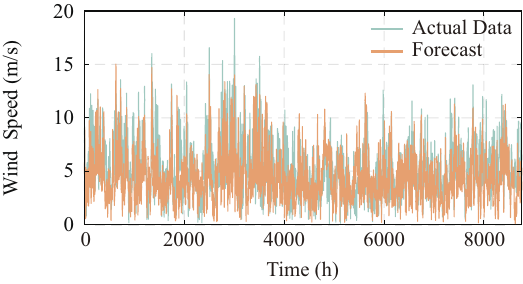}
		\vspace{-6pt}
		\caption{Wind speed}
		\label{fig:WindSpeed}
	\end{subfigure}
	%\qquad
	\begin{subfigure}{0.45\linewidth}
		\centering
		\includegraphics[width=0.95\linewidth]{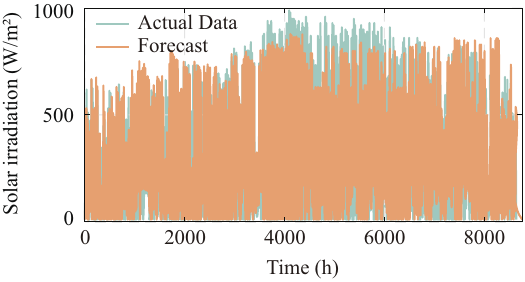}
		\vspace{-6pt}
		\caption{Solar irradiation}
		\label{fig:SolarIrra}
	\end{subfigure}
    \begin{subfigure}{0.45\linewidth}
		\centering
		\includegraphics[width=0.95\linewidth]{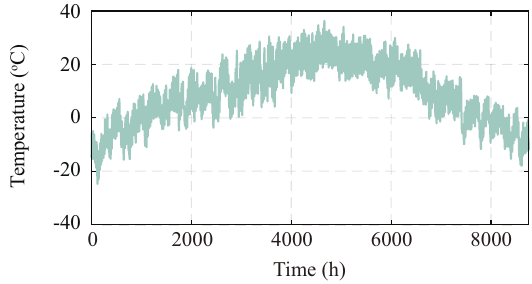}
		\vspace{-6pt}
		\caption{Ambient temperature}
		\label{fig:Temperature}
	\end{subfigure}
	\vspace{-6pt}
    \caption{Meteorological data collected from the site of the planned OReP2HS in the case study.}
	\label{fig:MeteorologicalData}
\end{figure}

The techno-economic parameters and base-case settings are presented in Tables~\ref{tab:TEparameter} and~\ref{tab:InitialData}. {\color{black}As no established operational standards currently exist for OReP2HSs, the frequency security limits in this case study are defined by referencing the permissible frequency ranges of large interconnected grids~\cite{FaVOR2010} and tailoring them based on the off-grid operating feature of the OReP2HS ~\cite{IEC2023,Farrokhabadi2019}. In addition, the ramping rate of each EL is conservatively set to 1\% of its rated load per second, which is lower than the experimentally demonstrated response capability of up to 10\% of the rated load per second~\cite{cheng2025power,sha2023low}. Therefore, the assumed ramping rate is well within practical engineering limits and can be achieved in real-world OReP2HS applications.} Fig. \ref{fig:MeteorologicalData} shows the meteorological data and on-site forecasts used in the simulations. Local ambient temperature, which affects the PV plant output, is obtained from the \emph{POWER Data Access Viewer} \cite{DataAccessViewer} with a time resolution of 1 hour. Observed and forecasted wind speeds and solar irradiance data with resolutions of 1 and 60 seconds, respectively, are sourced from Lingxiang Wind Farm and Jin'ao PV Plant in Inner Mongolia. {\color{black}These datasets serve as inputs to the simulations, which inherently account for the uncertainties of renewable power generation.} The 1-second resolution data are utilized to test grid-forming ability and ensure continuous operation during emergencies, while the 60-second resolution data are used for the annual energy balance test.
\begin{table}[tb]\footnotesize
  \renewcommand{\arraystretch}{1.15}
  \caption{Techno-economic parameters for LCOH evaluation}
  \label{tab:TEparameter}
  \vspace{-6pt}
  \centering
  \begin{tabular}{cccccc}
  \hline \hline
  {Facilities}  & {Number} & {Rating }  & {Investment cost} & {O$\&$M parameter $\lambda_j^\mathrm{O\&M}$}  & {Lifetime of facilities}  \\ %
  \hline

  {WT}  & {3} & {6,250 kW}  & \makecell[c]{5,000 CNY/kW} & \multirow{5}{*}{{\makecell[c]{2 \%}}}   & \multirow{5}{*}{{\makecell[c]{20 years}}}  \\

  {PV}  & {1} & {5,000 kW}  & \makecell[c]{4,000 CNY/kW} & \multirow{6}{*}{{\makecell[c]{}}}  &  \multirow{6}{*}{{\makecell[c]{}}} \\

  {AEL}  & {4} & {5,000 kW}  & \makecell[c]{3,500 CNY/kW} & \multirow{6}{*}{{\makecell[c]{}}} & \multirow{6}{*}{{\makecell[c]{}}}   \\

  {BESS}  & {1} & {To be evaluated}  & \makecell[c]{1,500 CNY/kWh} &\multirow{6}{*}{{\makecell[c]{}}} & \multirow{6}{*}{{\makecell[c]{}}}   \\

  {Transmission lines}  & {22 km} & \makecell[c]{35 kV}  & \makecell[c]{250,000 CNY/km} &\multirow{6}{*}{{\makecell[c]{}}} & \multirow{6}{*}{{\makecell[c]{}}}   \\

  \hline \hline
  \end{tabular}
\end{table}

\begin{table}[tb]\footnotesize
  \renewcommand{\arraystretch}{1.15}
  \caption{Base-case parameter settings for production simulation and BESS size search}
  \label{tab:InitialData}
  \vspace{-6pt}
  \centering
  \begin{tabular}{cc||cc}
  \hline \hline
  {Parameter}  & {Value} & {Parameter}  & {Value}    \\ %
  \hline
  {$T_{RS}$}    & {4 h}    & {$K_\mathrm{SOC}$}  & {0.0142}    \\ %
 % {$C_{\mathrm{H_{2}}}$}    & {30 CNY/kg}    & {$\beta$}  & {0.286}    \\ %
  {$C_{\mathrm{up,h}}^{\mathrm{EL}}$}  & {2 CNY}  & {$\beta$}  & {0.286}    \\ %
  {$C_{\mathrm{up,c}}^{\mathrm{EL}}$}  & {10 CNY} & {$c_{j}^{\mathrm{rep}}$}  & {900 CNY/kWh}    \\ %
  {$C_{\mathrm{down}}^{\mathrm{EL}}$}  & {5 CNY}& {$c_{j}^{\mathrm{rec}}$}  & {150 CNY/kWh}   \\ %
  {$s$}    & {1,000 CNY/MW}      & {$f_{\text{min}}^{\mathrm{PCC}}$}  & {46.5 Hz}   \\ %

  {$T_{\text{min}}^{\mathrm{down}}$}  & {1 h}    & $f_{\text{max}}^{\mathrm{PCC}}$  & {53.5 Hz}   \\ %

  {$k_{i}^{\mathrm{H_{2}}}$}    & {55.62 kWh/kg} & $V_{\text{min}}^{\mathrm{PCC}}$  & {31.5 kV}    \\ %

  {$SOC_{\text{min}}^{\mathrm{BES}}$}  & {10 \%}  & $V_{\text{max}}^{\mathrm{PCC}}$  & {38.5 kV}   \\ %
  {$SOC_{\text{max}}^{\mathrm{BES}}$}  & {90 \%}  &  $\Delta t_\mathrm{RS}$    & {5 min} \\
  $\Delta t_{\text{sim}\_\text{EB}}$      & {5$\times$10$^{-4}$ s}  &    $\Delta t_\mathrm{LF}$   &  {5 s} \\
  {$\Delta t_{\text{sim}\_\text{GF}}$}  & {4$\times$10$^{-5}$ s}  & \tabincell{c}{Ramping limit\\of the EL} & {0.05 MW/s} \\
  {$\Delta t_{\text{sim}\_N-1}$}  & {4$\times$10$^{-5}$ s}  &    &  \\
  \hline \hline
  \end{tabular}
\end{table}

\subsection{Base-Case Simulation and Evaluation LCOH Results}
\label{sec:SimResult}

\subsubsection{Production Simulation and BESS Sizing}
\label{sec:BESSsize}

\begin{figure}[H]
  \centering
  \includegraphics[width=3.25in]{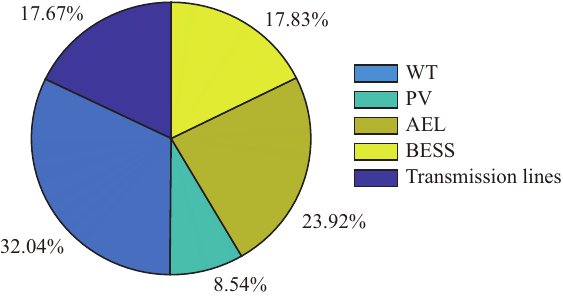}
  \vspace{-6pt}
  \caption{Proportion of investment and annual O\&M cost of devices in the OReP2HS under the optimal BESS size.}
  \label{fig:InvestProp}
\end{figure}

The optimal size of the grid-forming BESS is evaluated using the iterative search procedure described in Section \ref{sec:BatOptim} based on the proposed multi-timescale EMS shown in Fig. \ref{fig:EMS}. The sizing results and corresponding LCOH values are summarized in Table \ref{tab:SimResult-BESS}, indicating that the optimal BESS size is 6.8 MW/3.4 MWh.
The battery capacity of the BESS is equivalent to only 13.6\% of the maximal hourly total energy output of the WTs and PV plant. The annual costs, including the annualized capital cost and annual O\&M cost of the BESS, account for 17.83\% of the total annual cost, as illustrated in Fig. \ref{fig:InvestProp}. {\color{black}Moreover, Table~\ref{tab:cost_shares} indicates that the O\&M cost of the BESS contributes more significantly to the overall LCOH than its investment cost.}
\begin{table}[tb]\footnotesize
  \renewcommand{\arraystretch}{1.15}
  \caption{Base-case result of BESS sizing and LCOH evaluation}
  \label{tab:SimResult-BESS}
  \vspace{-6pt}
  \centering
  \begin{tabular}{cccc}
  \hline \hline
   \tabincell{c}{Battery size} & \tabincell{c}{Charging/Discharging rate}  & \tabincell{c}{Yearly degradationy}   & \tabincell{c}{LCOH (CNY/kg)} \\ %
  \hline
   \makecell[c]{3.4 MWh} & \makecell[c]{2C}  & \makecell[c]{4.87 \%} & \makecell[c]{33.212} \\ %
  \hline \hline
  \end{tabular}
\end{table}

\begin{table}[H]
  \centering
  \footnotesize
  \caption{\color{black}The contribution of each device’s investment and O\&M costs to the LCOH}
  \label{tab:cost_shares}
    \vspace{-6pt}
  {\color{black}
  \renewcommand{\arraystretch}{1.15}
  \begin{tabular}{c*{10}{c}} % <-- 将 l 改为 c
    \hline\hline
    \multirow{2}{*}{Share} &
    \multicolumn{2}{c}{BESS} &
    \multicolumn{2}{c}{WT} &
    \multicolumn{2}{c}{PV} &
    \multicolumn{2}{c}{AEL} &
    \multicolumn{2}{c}{Transmission lines} \\
    \cline{2-11}
    & Invest. & O\&M & Invest. & O\&M & Invest. & O\&M & Invest. & O\&M & Invest. & O\&M \\
    \hline
    Percentage (\%) &
    3.20 & 14.63 &
    26.02 & 6.02 &
    6.93 & 1.61 &
    19.42 & 4.50 &
    13.35 & 3.32 \\
    \hline\hline
  \end{tabular}}
\end{table}

\begin{figure}[tb]
  \centering
    \includegraphics[width=0.95\textwidth]{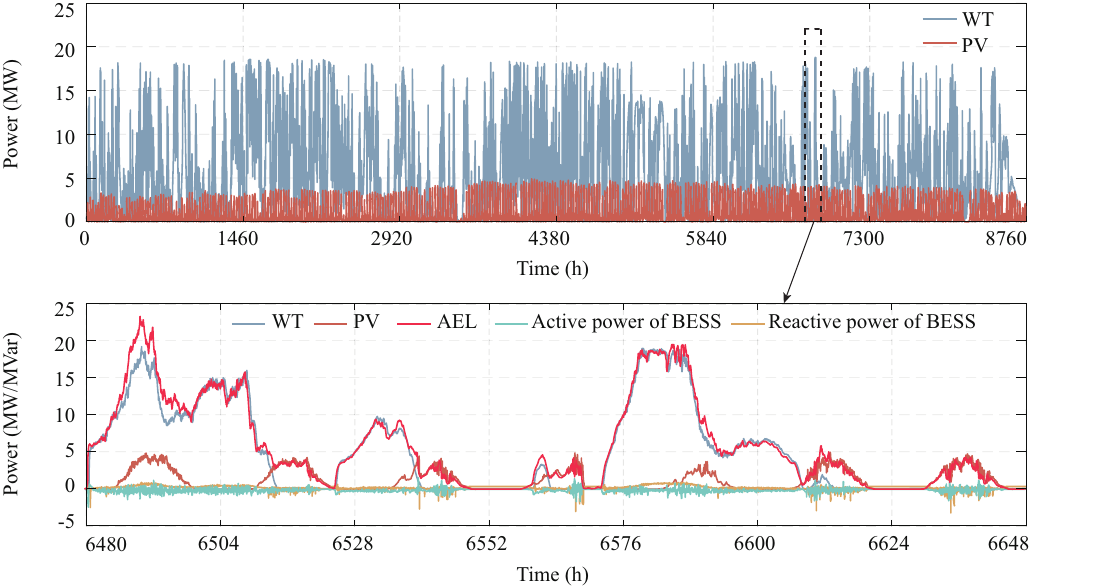}
    \vspace{-6pt}
    \caption{Power of PV, WTs, BESS and AELs in the year-long production simulation.}
    \label{fig:8760hPower}
  \end{figure}

Fig. \ref{fig:8760hPower} displays the production simulation results over a reference year, with detailed results from October 1st to October 7th. This period illustrates various operating conditions, including complementary operation, simultaneous fluctuations, and outages of the WT and PV outputs.

First, we discuss the performance of the proposed EMS. With the RS and LF strategies, the electrolytic load tracks the fluctuations in renewable power, optimizing the use of renewable energy and minimizing the required BESS capacity.
For instance, during hours 6,596 to 6,606, only the wind turbines (WTs) are operational, and the AELs closely follow wind power fluctuations. However, from hours 6,608 to 6,619, the AELs struggle to match the rapidly changing photovoltaic (PV) output due to its steep ramping rate. Under these conditions, the BESS must provide sufficient transient power support under grid-forming control to maintain power balance. In such scenarios, the commonly used minute-level resolution models in existing studies, {\color{black}which assume that power balance can be maintained within each time step, become invalid for PV-dominated OReP2HSs.} This leads to significant deviations in the estimated BESS sizing from actual system needs.

\begin{figure}[htbp]
  \centering
  \begin{subfigure}{\textwidth}
    \centering
    \includegraphics[width=0.95\textwidth]{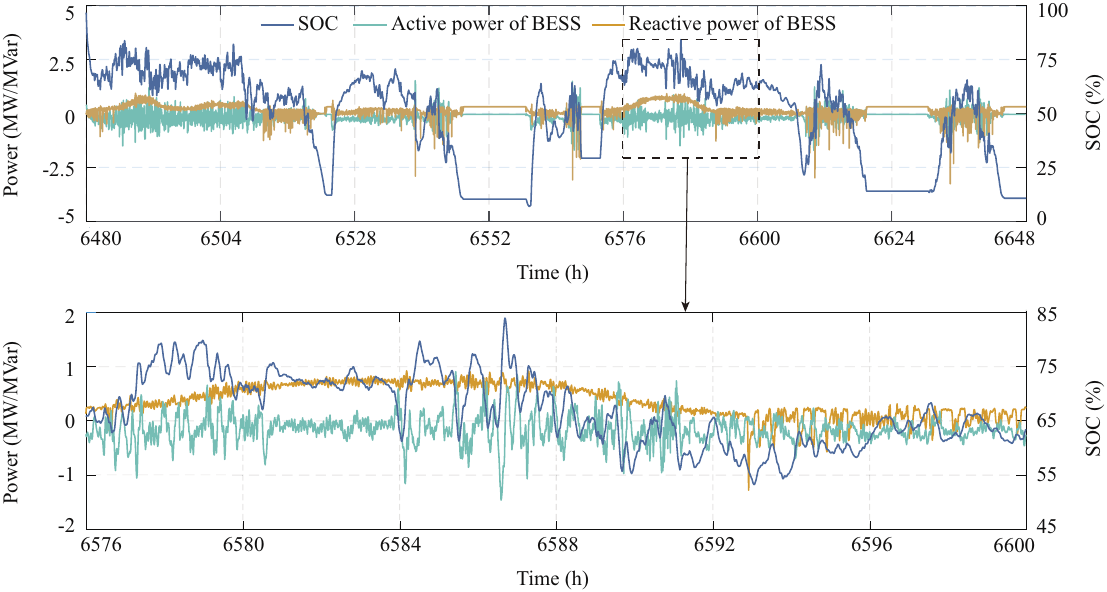}
    \vspace{-6pt}
    \caption{The SOC and active/reactive power of the grid-forming BESS in the year-long production simulation.}
    \label{fig:BatSOC}
  \end{subfigure}
  \vspace{3pt}

  \begin{subfigure}{\textwidth}
    \centering
    \includegraphics[width=0.95\textwidth]{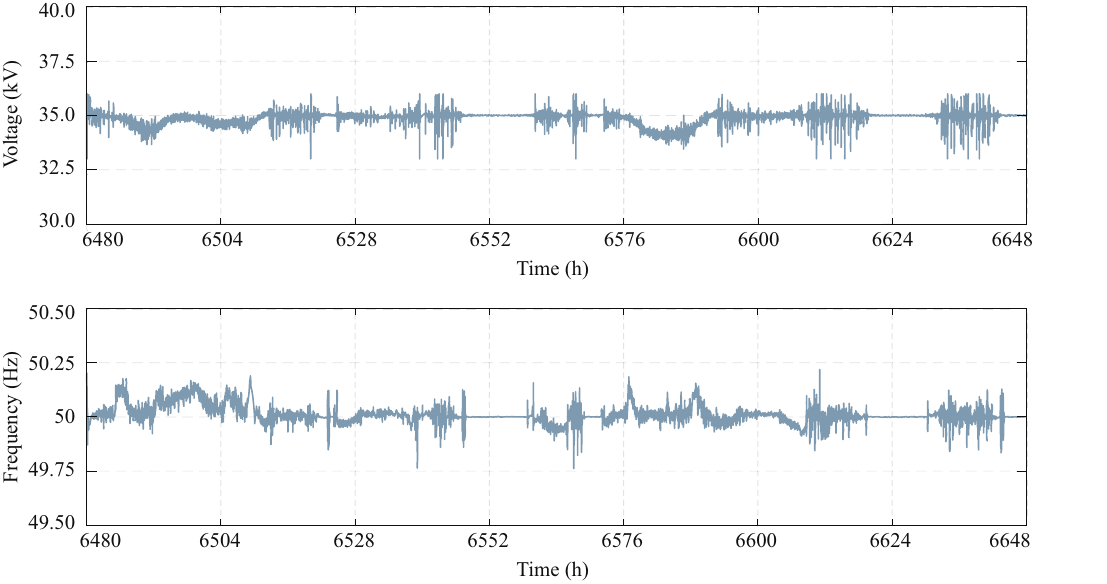}
    \vspace{-6pt}
    \caption{Frequency and voltage at PCC from October 1st to 7th in the production simulation.}
    \label{fig:VoltageAndFrequency}
  \end{subfigure}

  \caption{Detailed production simulation results under the proposed EMS.}
  \label{fig:ProductionSimulationAll}
\end{figure}

The SOC and power of the BESS, as shown in Fig. \ref{fig:BatSOC}, further support this observation. The BESS power fluctuates rapidly to smooth out high-frequency renewable power variations occurring within the transient timescale that existing models fail to capture. Notably, active power exhibits more frequent fluctuations than reactive power, highlighting the importance of frequency security in OReP2HSs. Therefore, a specific EH strategy has been designed to ensure frequency security under emergency conditions (see Section \ref{sec:SimCODE}). Meanwhile, the SOC remains within the acceptable range (10\%–90\%), demonstrating the effectiveness of the SOC correction strategy embedded in the LF control. The SOC trajectory mirrors the profile of renewable generation under our proposed EMS, indicating a temporal correlation discussed further in Section \ref{sec:Correlation}.
%\begin{figure}[tb]
%  \centering
%  \includegraphics[width=5.9in]{8760hPower.eps}
%  \vspace{-6pt}
%  \caption{Power of PV, WTs, BESS and AELs in the year-long production simulation.}
%  \label{fig:8760hPower}
%\end{figure}
%
%\begin{figure}[h]
%  \centering
%  \includegraphics[width=5.9in]{BatSOC.eps}
%  \vspace{-6pt}
%  \caption{The SOC and active/reactive power of the grid-forming BESS in the year-long production simulation.}
%  \label{fig:BatSOC}
%\end{figure}
%
%\begin{figure}[h]
%  \centering
%  \includegraphics[width=5.9in]{VoltageAndFrequency.eps}
%    \vspace{-6pt}
%  \caption{Frequency and voltage at PCC from October 1st to 7th in the production simulation.}
%  \label{fig:VoltageAndFrequency}
%\end{figure}
Finally, the frequency and voltage at PCC from October 1st to 7th are presented in Fig. \ref{fig:VoltageAndFrequency}, which aligns with the active and reactive power patterns of the BESS in Fig. \ref{fig:BatSOC}. This demonstrates that frequency and voltage dynamics in OReP2HSs are strongly coupled with the active and reactive power responses of the BESS. Consequently, BESS sizing is not only a matter of energy and power balancing across timescales but also critically influences system voltage and frequency dynamic. Voltage and frequency constraints must therefore be considered during BESS configuration. Table \ref{tab:SimResult-VF} confirms that frequency and voltage remain within acceptable limits throughout the period, validating the effectiveness of the proposed EMS and the adequacy of the optimized BESS sizing in supporting system stability.

\begin{table}[h]\footnotesize
  \renewcommand{\arraystretch}{1.15}
  \caption{Technical indexes of the frequency and voltage at PCC}
  \label{tab:SimResult-VF}
  \vspace{-6pt}
  \centering
  \begin{tabular}{ccccccc}
  \hline \hline
  \makecell[c]{Index}  & \makecell[c]{Rated value}  & \makecell[c]{Allowable interval} & \makecell[c]{Peak}  & \makecell[c]{Nadir} &  \makecell[c]{Average}  & \tabincell{c}{Ratio of the maximal \\deviation to the rating} \\ %
  \hline
  \makecell[c]{Frequency}  & \makecell[c]{50 Hz}  & \makecell[c]{(46.5 Hz, 53.5 Hz)}   & \makecell[c]{50.22 Hz} &  \makecell[c]{49.77 Hz}  & \makecell[c]{50.02 Hz} & \makecell[c]{0.46\%} \\ %
  \makecell[c]{Voltage}  & \makecell[c]{35 kV}  & \makecell[c]{(31.5 kV, 38.5 kV)}   & \makecell[c]{36.00 kV} &  \makecell[c]{33.01 kV}  & \makecell[c]{34.89 kV} & \makecell[c]{5.69\%} \\ %
  \hline \hline
  \end{tabular}
\end{table}

\subsubsection{Operation During Emergencies}
\label{sec:SimCODE}
In this section, we conduct a WT tripping test to analyze the performance of the proposed EH strategy. The BESS size is set at 6.8 MW/3.4 MWh, and the following settings are compared:

\textbf{Setting A:} The EH is automatically activated based on its activation criteria listed in Table \ref{tab:Emergency}, as in the proposed EMS used in search for the optimal BESS size in Section \ref{sec:BESSsize}.

\textbf{Setting B:} The EH is removed, and other sub-strategies of the EMS remain consistent with Setting A.
\begin{figure}[tbh]
  \centering
  \includegraphics[width=5.7in]{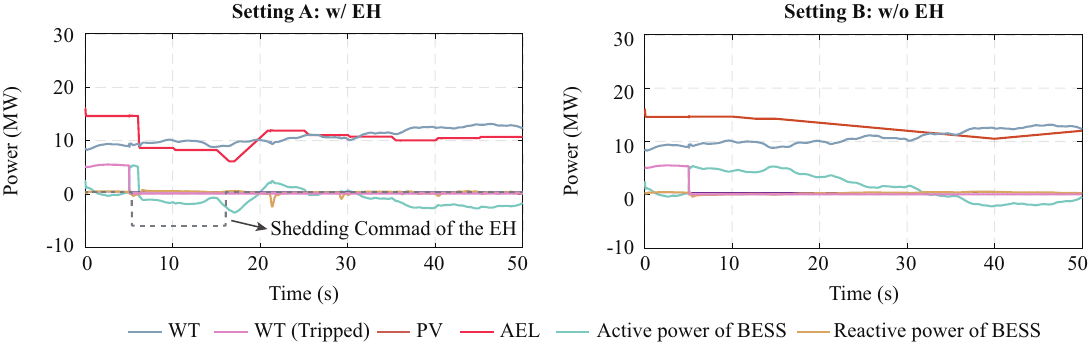}
  \vspace{-6pt}
  \caption{The power variations of PV, WTs, BESS and AELs when a WT is tripped.}
  \label{fig:CODE_Power}
\end{figure}
\begin{figure}[tbh]
  \centering
  \includegraphics[width=5.7in]{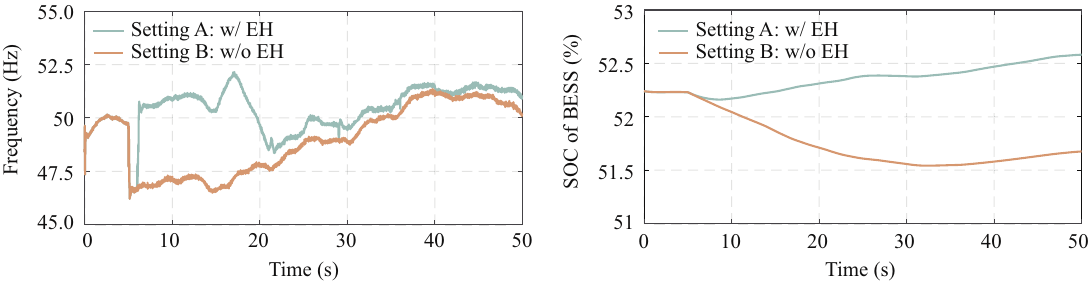}
  \vspace{-6pt}
  \caption{The frequency and SOC variation when a WT tripped.}
  \label{fig:CODE_SOC-f}
\end{figure}
The simulation results with a WT tripping at 5 s are depicted in Figs. \ref{fig:CODE_Power} and \ref{fig:CODE_SOC-f}. When the WT trips at 5 s, the frequency drops sharply, while the BESS output increases rapidly to prevent frequency collapse. The EH is activated when the frequency deviation reaches 0.35 Hz. Subsequently, power shedding commands are sent to all AELs. The electrolytic load decreases from 15 MW to 10 MW with a maximal ramping rate of 0.5 MW/s after 1 s, allowing the frequency to recover to its rated value within 1.5 s. Meanwhile, the BESS output decreases back to the initial state.

By contrast, without the EH, AELs cannot be regulated to respond promptly to the power loss. Their loads slowly reduce from 15 MW to 10 MW as commanded by the LF strategy. This process takes nearly 29 s, causing the slow recovery of frequency. During this period, the BESS maintains its output to balance the power loss, leading to a decrease of 0.8\% in its SOC.

Note that the system remains stable when the WT is tripped under both settings (with and without the EH strategy). The proposed EH strategy leverages the flexibility of AELs to mitigate SOC variations during emergencies, thereby reducing the required BESS capacity.

{\color{black}We also want to remark that the proposed lookup-table-based EH is designed for practical and straightforward implementation in real-world projects. Equipping the system with more advanced emergency handling mechanisms (e.g., adaptive strategies) may further improve performance and reduce the required BESS capacity. Although the proposed EH strategy is reasonably effective, further exploration is needed to fully understand the potential and limitations of the OReP2HS in addressing emergencies.}

\subsubsection{Correlation Between the SOC and Renewable Outputs}
\label{sec:Correlation}
\begin{figure}[tb]
  \centering
  \includegraphics[width=5.9in]{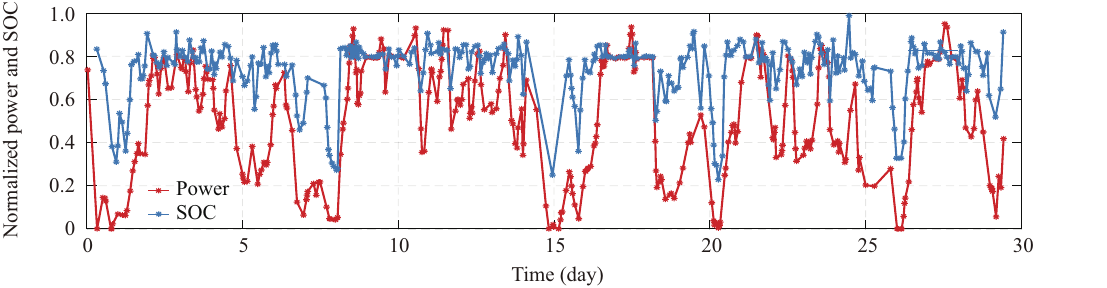}
  \vspace{-6pt}
  \caption{Temporal variations of local extrema in the SOC of the BESS and renewable output in March.}
  \label{fig:isolaMAX}
\end{figure}

As observed in Fig. \ref{fig:BatSOC}, under the proposed EMS, the SOC trajectory of the BESS closely follows the fluctuations in renewable power output. Fig. \ref{fig:isolaMAX} further illustrates the local extrema in SOC and renewable power for March. Whenever renewable power reaches a local extremum, a corresponding peak in SOC occurs shortly afterward, indicating a temporal correlation between their fluctuations. This temporal similarity motivates an investigation into the correlation between renewable power and the grid-forming BESS behavior, with the aim of exploring potential strategies to reduce the BESS size.

\begin{figure}[htb]
	\centering
	\begin{subfigure}{0.45\linewidth}
		\centering
		\includegraphics[width=0.90\linewidth]{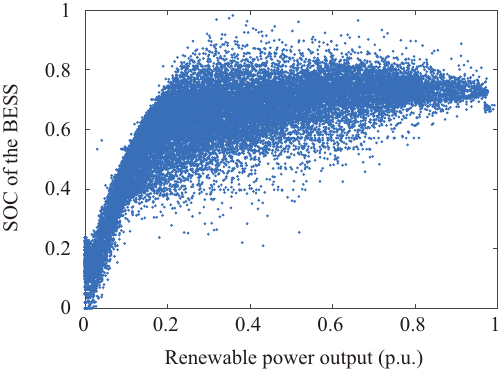}
		\vspace{-6pt}
		\caption{Scatter}
		\label{fig:Scatter}
	\end{subfigure}
    %\qquad
	\begin{subfigure}{0.5\linewidth}
		\centering
		\includegraphics[width=0.94\linewidth]{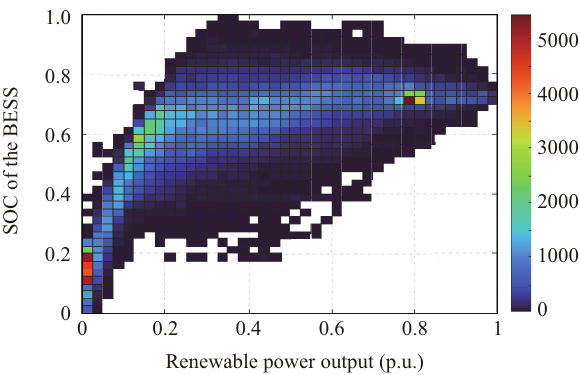}
		\vspace{-6pt}
		\caption{Local statistical results}
		\label{fig:Stastical}
	\end{subfigure}
	\vspace{-6pt}
    \caption{Statistical relationship between the SOC of the BESS and renewable power.}
	\label{fig:StatisticalResults}
\end{figure}

{\color{black}A statistical analysis of the yearly (8760 hours) SOC time series $\{\boldsymbol{\xi}_{\mathrm{SOC}}(t)\}$ and renewable power output time series $\{\boldsymbol{\xi}_{\mathrm{re}}(t)\}$ with 1-min time resolution is presented in Fig. \ref{fig:StatisticalResults}.} As shown in Fig. \ref{fig:Scatter}, at the micro level, the SOC exhibits a positive and nonlinear relationship with renewable power. The SOC changes rapidly in the lower renewable output range (normalized power ranging from 0 to 0.4). During this interval, the precision of fast forecasting integrated into the LF strategy influences the outputs and SOC of the BESS. Frequent and irregular ramps of the renewable power complicate accurate forecasting, necessitating frequent BESS output adjustments, as shown in the local statistics in Fig. \ref{fig:Stastical}. This phenomenon leads to an almost linear relationship between SOC and renewable power in this range.

By contrast, when renewable generation operates at a higher output interval (normalized power from 0.4 to 1), the outputs are relatively stable. During these periods, AELs consume nearly all renewable outputs, leaving the BESS to balance minor power fluctuations, and causing the SOC of the BESS to become less correlated with renewable power.

Additionally, Fig. \ref{fig:Stastical} shows that the normalized SOC distribution is mainly concentrated around 0.15 and 0.75, indicating that the BESS frequently operates in deep charge or discharge states. This explains the rapid BESS degradation rate of 4.87\% annually as reported in Section \ref{sec:BESSsize}.

{\color{black}Furthermore, we employ the Spearman rank correlation coefficient $r_s(\boldsymbol{\xi}_{\mathrm{SOC}}, \boldsymbol{\xi}_{\mathrm{re}})$ to quantify the monotonic dependence between the SOC and renewable generation, as defined by:
\begin{equation}
r_s(\boldsymbol{\xi}_{\mathrm{SOC}}, \boldsymbol{\xi}_{\mathrm{re}}) =
\frac{\sum_{i=1}^{n}(R_i-\bar{R})(S_i-\bar{S})}
{\sqrt{\sum_{i=1}^{n}(R_i-\bar{R})^2}\sqrt{\sum_{i=1}^{n}(S_i-\bar{S})^2}}, \label{rannk}
\end{equation}
where $R_i$ and $S_i$ are the rank values of $\{\boldsymbol{\xi}_{\mathrm{SOC}}(t)\}$ and $\{\boldsymbol{\xi}_{\mathrm{re}}(t)\}$ at the $i$th time step, respectively, and $\bar{R}$ and $\bar{S}$ denote their mean ranks.}

The calculation results listed in Table~\ref{tab:Rank} suggest that, at the macro level, a notable temporal correlation exists between the BESS SOC and renewable generation, SOC prediction based on renewable power forecasts can enable actions to mitigate battery degradation, optimize BESS energy management, and potentially further reduce the BESS size.

\begin{table}[tb]\footnotesize
  \renewcommand{\arraystretch}{1.15}
  \caption{Rank correlation between the SOC of BESS and renewable power in each month}
  \vspace{-6pt}
  \label{tab:Rank}
  \centering
  \begin{tabular}{ccccccccccccc}
  \hline \hline
  {}  & {Mar.}  & {Apr.} & {May}  & {Jun.} &  {Jul.}  & {Aug.} & {Sep.}  & {Oct.} & {Nov.}  & {Dec.} & Jan. & {Feb.}\\ %
  \hline
  Rank correlation  & 0.855   & 0.851   & 0.871  &  0.902  & 0.870 & 0.835  & 0.833 & 0.858  & 0.845   & 0.896  &  0.858  & 0.847 \\ %
  \hline \hline
  \end{tabular}
\end{table}

\subsection{Comparison of the Optimal BESS Size with Different EMS Under Stability Constraints}
\label{sec:EMSCompair}
As summarized in Section \ref{sec:review}, existing EMSs supporting the design of OReP2HS are predominantly based on rule-based or optimization-based approaches with time resolutions ranging from at least 5 minutes to 1 hour. The components of the OReP2HS are optimized based on these EMSs considering only the power balance of PV, WTs, AELs, and BESS at discrete time steps (every 5 minutes or hourly), without accounting for the intra-step power balance. However, BESS may overcharge or over-discharge between these discrete time steps, losing its grid-forming ability and causing instability. Thus, the design of OReP2HS based on EMS with 5-min or hourly resolutions in the may lead to the underestimation of the BESS size required to ensure the reliable operation of OReP2HS.

To this end, this section compares the optimal battery capacities reported in the state-of-the-art studies \cite{ibanez2022simulation,abdin2024feasibility} under EMS strategies Benchmark 1 and 2 with a time resolution of 5 minutes, and conducts two tests: 1) using the proposed high-fidelity model (with a time resolution of 0.04 ms) to test whether the battery capacities reported in existing literature can support the stable operation of an OReP2H system over a full year (8760 hours) simulation; 2) if frequency or voltage instability occurs during the 8760-hour simulation using the reported capacity, applying the proposed optimization method (as shown in Fig. \ref{fig:Solving}, which considers the voltage and frequency security as constraints) to re-evaluate the minimum battery capacity required for stable operation under the literature's EMS logic.

\textbf{Benchmark 1:} A rule-based EMS from \cite{ibanez2022simulation,ibanez2023off}. The strategy is as follows: if the renewable power exceeds the minimal load of AELs, the AELs are started up, and then the electrolytic load tracks the renewable power with a time step of 5 minutes. When the maximal loads of AELs are reached, the rest of the power charges the BESS. When the combined power from PV, WTs, and the BESS is insufficient to keep the AELs at their minimum load, the AELs are shut down.

\textbf{Benchmark 2:} An MILP-based EMS similar to those in \cite{wang2023optimising,abdin2024feasibility}. AELs and the BESS track the renewable power only follow an MILP-based scheduling results. For fair comparison, we take (\ref{eq-1})--(\ref{eq-17}) introduced in Section \ref{sec:RS} as the MILP model and set the time step as 5 minutes.

{\color{black}According to the literature review in Section~\ref{sec:review}, the rule-based and MILP-based EMSs represent the two most widely adopted EMS for designing OReP2HS. Therefore, they are selected here as representative benchmarks. To ensure fairness, both benchmark methods are configured with the highest temporal resolution (5 minutes) reported in the state-of-the-art studies. However, this resolution remains insufficient to capture the transient behaviors of OReP2HS, making it impossible to incorporate frequency and voltage security constraints during BESS sizing.}

\begin{figure}[tb]
    \centering
    \begin{subfigure}[t]{0.45\linewidth}
        \includegraphics[width=0.90\linewidth]{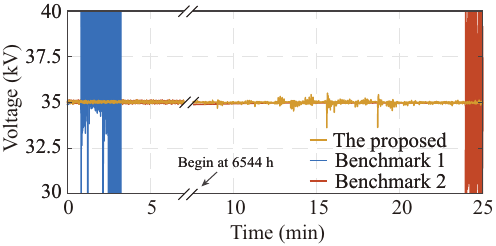}
        \vspace{-6pt}
        \caption{Voltage at PCC}
        \label{fig:VotgCompare}
    \end{subfigure}
    \hspace{0.75em}
    \begin{subfigure}[t]{0.45\linewidth}
        \includegraphics[width=0.90\linewidth]{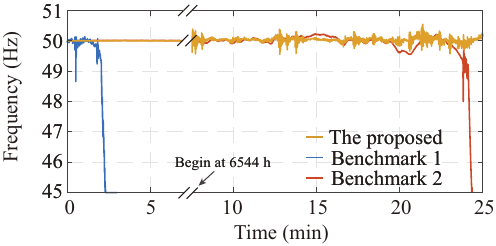}
        \vspace{-6pt}
        \caption{Frequency at PCC}
        \label{fig:FreqCompare}
    \end{subfigure}
    \vspace{-6pt}
    \caption{Comparison of the simulation results under different EMS.}
    \label{fig:EMSComparison-1}
\end{figure}

\begin{figure}[htb]
    \centering
    \begin{subfigure}[t]{0.45\linewidth}
        \includegraphics[width=0.90\linewidth]{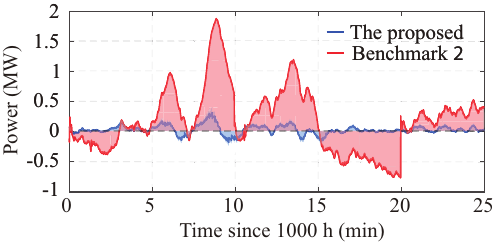}
        \vspace{-6pt}
        \caption{Active power of BES}
        \label{fig:BesOutCompare}
    \end{subfigure}
    \hspace{0.75em}
    \begin{subfigure}[t]{0.45\linewidth}
        \includegraphics[width=0.90\linewidth]{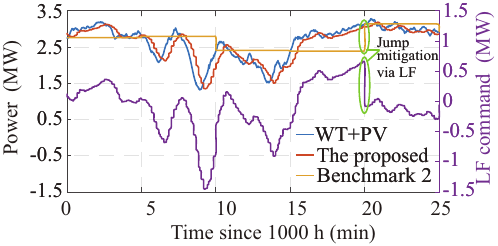}
        \vspace{-6pt}
        \caption{Renewable power and load command of AEL}
        \label{fig:LoadPowerDecomposed}
    \end{subfigure}
    \vspace{-6pt}
    \caption{Comparison of the simulation results under Benchmark 2 and the proposed methods.}
    \label{fig:EMSComparison-2}
\end{figure}

%\textbf{The proposed :} the production simulation of OReP2HS is performed based on \textit{the proposed EMS}.
The simulation results are summarized in Table \ref{tab:SimResultEMS-BESS}. For all EMS strategies reported in the literature, frequency and voltage instabilities were observed during the 8760-hour simulation (see Fig. \ref{fig:EMSComparison-1}), primarily due to the lower time resolution of the optimal under Benchmark 1 and 2. The components of the OReP2HS are optimized based on these EMSs considering only the power balance of PV, WTs, AELs, and BESS at a 5-min discrete time steps, without accounting for the intra-step power balance. However, in practice, renewable generation fluctuations within these intervals result in power imbalances that must be mitigated through BESS charge/discharge operations. When evaluated using the proposed high-fidelity simulation model, these intra-step power fluctuations lead BESS overcharge or over-discharge, losing its grid-forming ability and causing instability. After re-evaluation using the proposed method, the battery capacities required to satisfy stability constraints are higher than those originally reported, indicating that prior studies underestimated the minimum BESS size.
\begin{table}[tb]\footnotesize
	\renewcommand{\arraystretch}{1.1}
	\caption{Optimal BESS size and LCOH under different EMS}
	\vspace{-6pt}
	\label{tab:SimResultEMS-BESS}
	\centering
	\begin{tabular}{cccccc}
		\hline \hline
		\multirow{2}*{\makecell[c]{EMS}}
		& \multicolumn{2}{c}{The optimal configuration ratio of battery ${\rho_{\text{BES}}}^{1} $}
		&  \multirow{2}*{\tabincell{c}{Charging/\\Discharging rate}}
		&  \multirow{2}*{\tabincell{c}{Yearly battery\\ degradation}}
		&  \multirow{2}*{\tabincell{c}{LCOH \\(CNY/kg)}} \\  \cline{2-3}%\cline{2-4}
		&  \tabincell{c}{Reported in literature \\ (Stability maintained$^{2}$?) }
		&  \tabincell{c}{Re-evaluated in this paper \\ (Stability maintained$^{2}$?) } &  &   \\ %		
       \hline
		Benchmark 1 & \tabincell{c}{3.6\% (0.9 MWh) \cite{ibanez2022simulation} \\ (\text{\sffamily X})} & \tabincell{c}{69.44\% (17.36 MWh)\\ (\checkmark)} & \makecell[c]{1C}
        & 3.91\% & 56.958 \\%
		{\makecell[c]{Benchmark 2}} & \tabincell{c}{ 26.28\%$\sim$40.11\% \\(6.67$\sim$10.03 MWh) \cite{abdin2024feasibility} \\ (\text{\sffamily X})} &{\makecell[c]{55.67\% (13.94 MWh) \\ (\checkmark)}}
        & {\makecell[c]{1C}}  &{\makecell[c]{3.55 \%}} &{\makecell[c]{51.737} }  \\ %
		\bf The proposed &  /  & \tabincell{c}{13.6\% (3.4 MWh)  \\ (\checkmark)} & 2C &  4.87 \% & 33.212 \\ %
		\hline \hline
         \multicolumn{6}{l}{$^{1}$ $\rho_{\text{BES}}=S_{\text{BES}}/(S_{\text{PV}}+\sum_{i=1}^{m}S_{\text{WT},i})$}\\
        \multicolumn{6}{l}{$^{2}$ Testing system frequency and voltage stability over an 8760-hour operation using the proposed high-fidelity simulation}\\
       \multicolumn{6}{l}{ model under different battery capacity.}

	\end{tabular}
\end{table}
Moreover, the re-evaluated battery size and LCOH obtained using the proposed EMS are significantly lower than those of Benchmark 1 and 2, demonstrating superior economic performance. This improvement is mainly attributed to the proposed LF controller, which integrates a fast ARMA process-based renewable power prediction module and enables 5-second high-frequency adjustments to correct the baseline power determined by the RS. This coordination ensures that the AELs more closely follow the actual renewable generation outputs, thereby reducing the system’s reliance on batteries to compensate for source-load power mismatch. As shown in Fig. \ref{fig:BesOutCompare}, the proposed method results in lower battery charge/discharge power fluctuations and reduced total energy consumption compared to Benchmark 2.

Besides, Fig. \ref{fig:LoadPowerDecomposed} indicates that the RS and LF modules operate in a well-coordinated manner, and the load power of AEL exhibits smooth transitions without abrupt changes caused by signal mismatch between the RS and LF.

\subsection{Sensitivity Analysis}
\label{sec:Sensitivity}
Taking advantage of the load flexibility of AELs, part of the energy balance regulation requirements for the grid-forming BESS is replaced, thus reducing the BESS size and LCOH. The load flexibility of AELs is significantly influenced by two key parameters in the EMS: the time step of the LF, $\Delta t_\mathrm{LF}$, and the load ramping limit of AELs. Therefore, a sensitivity analysis of these two parameters is performed to provide a comprehensive understanding of the limitations of the proposed EMS and to explore the potential for further reducing the BESS size and LCOH.

For the sensitivity analysis, six parameter setting groups are considered for comparison. In each group, the ramping limit of AELs is kept constant while the time step of the LF is increased stepwise. The results of the sensitivity analysis are listed in Table \ref{tab:Sensitivity}.

\begin{table}[tbph]\footnotesize
  \renewcommand{\arraystretch}{1.08}
  \caption{Sensitivity analysis result of the optimal BESS size and LCOH with different parameter settings}
  \vspace{-6pt}
  \label{tab:Sensitivity}
  \centering
  \begin{tabular}{cccccc}
  \hline \hline
  \tabincell{c}{Parameter \\setting} & \tabincell{c}{Time step of \\ LF $\tau$}  & \tabincell{c}{Ramp limit \\of AEL (MW/s)}  & \tabincell{c}{Optimal battery size (MWh)/ \\ Charging and discharging rate}  & \tabincell{c}{Yearly degradation \\ of battery}   & \tabincell{c}{LCOH \\(CNY/kg)} \\ %
  \hline
  \textbf{Base}  & \textbf{5 s}  &\textbf{0.05}  & \textbf{3.40/2C}  & \textbf{4.87\%}  & \textbf{{33.212} } \\ %
  \hline
  \multirow{5}{*}{{\makecell[c]{\uppercase\expandafter{\romannumeral 1}}}} & {10 s}  &{0.05}  & {3.75/2C}  & {4.88\%}  & {33.767}\\
  \multirow{5}{*}{{\makecell[c]{}}} & {15 s}  &{0.05}  & {4.04/2C}  & {4.91\%}  & {34.281}\\
  \multirow{5}{*}{{\makecell[c]{}}} & {30 s}  &{0.05}  & {4.62/2C}  & {5.01\%}  & {38.063}\\
  \multirow{5}{*}{{\makecell[c]{}}} & {60 s}  &{0.05}  & {5.07/2C}  & {4.97\%}   & {39.951}\\
  \multirow{5}{*}{{\makecell[c]{}}} & {90 s}  &{0.05}  & {5.73/2C}  & {4.99\%}  & {42.211}\\
  \hline

  \multirow{6}{*}{{\makecell[c]{\uppercase\expandafter{\romannumeral 2}}}} & {5 s}  &{0.1}  & {3.05/2C}  & {4.94\%}   & {31.551}\\
  \multirow{6}{*}{{\makecell[c]{}}} & {10 s}  &{0.1}  & {3.18/2C}  & {4.94\%}   & {32.777}\\
  \multirow{6}{*}{{\makecell[c]{}}} & {15 s}  &{0.1}  & {3.56/2C}  & {5.04\%}   & {33.437}\\
  \multirow{6}{*}{{\makecell[c]{}}} & {30 s}  &{0.1}  & {3.96/2C}  & {5.04\%}  & {35.133}\\
  \multirow{6}{*}{{\makecell[c]{}}} & {60 s}  &{0.1}  & {4.14/2C}  & {5.05\%}   & {35.891}\\
  \multirow{6}{*}{{\makecell[c]{}}} & {90 s}  &{0.1}  & {4.87/2C}  & {5.12\%}   & {38.715}\\
  \hline

  \multirow{6}{*}{{\makecell[c]{\uppercase\expandafter{\romannumeral 3}}}} & {5 s}  &{0.2}  & {2.78/2C}  & {4.94\%}   & {31.018}\\
  \multirow{6}{*}{{\makecell[c]{}}} & {10 s}  &{0.2}  & {2.89/2C}  & {4.93\%}   & {31.275}\\
  \multirow{6}{*}{{\makecell[c]{}}} & {15 s}  &{0.2}  & {3.24/2C}  & {5.03\%}   & {32.881}\\
  \multirow{6}{*}{{\makecell[c]{}}} & {30 s}  &{0.2}  & {3.60/2C}  & {5.11\%}   & {33.507}\\
  \multirow{6}{*}{{\makecell[c]{}}} & {60 s}  &{0.2}  & {3.76/2C}  & {5.12\%}   & {34.785}\\
  \multirow{6}{*}{{\makecell[c]{}}} & {90 s}  &{0.2}  & {4.32/2C}  & {5.14\%}   & {36.559}\\
  \hline

  \multirow{6}{*}{{\makecell[c]{\uppercase\expandafter{\romannumeral 4}}}} & {5 s}  &{0.3}  & {2.66/2C}  & {4.92\%}  & {30.872}\\
  \multirow{6}{*}{{\makecell[c]{}}} & {10 s}  &{0.3}  & {2.80/2C}  & {4.97\%}   & {31.116}\\
  \multirow{6}{*}{{\makecell[c]{}}} & {15 s}  &{0.3}  & {3.14/2C}  & {5.05\%}  & {32.707}\\
  \multirow{6}{*}{{\makecell[c]{}}} & {30 s}  &{0.3}  & {3.49/2C}  & {5.13\%}   & {33.316}\\
  \multirow{6}{*}{{\makecell[c]{}}} & {60 s}  &{0.3}  & {3.64/2C}  & {5.14\%}  & {33.977}\\
  \multirow{6}{*}{{\makecell[c]{}}} & {90 s}  &{0.3}  & {4.15/2C}  & {5.14\%}   & {35.463}\\
  \hline

  \multirow{6}{*}{{\makecell[c]{\uppercase\expandafter{\romannumeral 5}}}} & {5 s}  &{0.4}  & {2.66/2C}  & {4.97\%}   & {30.799}\\
  \multirow{6}{*}{{\makecell[c]{}}} & {10 s}  &{0.4}  & {2.80/2C}  & {4.97\%}   & {31.115}\\
  \multirow{6}{*}{{\makecell[c]{}}} & {15 s}  &{0.4}  & {3.10/2C}  & {5.10\%}   & {32.634}\\
  \multirow{6}{*}{{\makecell[c]{}}} & {30 s}  &{0.4}  & {3.43/2C}  & {5.11\%}   & {33.211}\\
  \multirow{6}{*}{{\makecell[c]{}}} & {60 s}  &{0.4}  & {3.59/2C}  & {5.15\%}  & {33.699}\\
  \multirow{6}{*}{{\makecell[c]{}}} & {90 s}  &{0.4}  & {4.07/2C}  & {5.16\%}   & {35.324}\\
  \hline

  \multirow{6}{*}{{\makecell[c]{\uppercase\expandafter{\romannumeral 6}}}} & {5 s}  &{0.5}  & {2.66/2C}  & {4.97\%}   & {30.799}\\
  \multirow{6}{*}{{\makecell[c]{}}} & {10 s}  &{0.5}  & {2.80/2C}  & {4.97\%}   & {31.109}\\
  \multirow{6}{*}{{\makecell[c]{}}} & {15 s}  &{0.5}  & {3.08/2C}  & {5.10\%}   & {32.602}\\
  \multirow{6}{*}{{\makecell[c]{}}} & {30 s}  &{0.5}  & {3.39/2C}  & {5.11\%}  & {33.142}\\
  \multirow{6}{*}{{\makecell[c]{}}} & {60 s}  &{0.5}  & {3.56/2C}  & {5.15\%}   & {33.437}\\
  \multirow{6}{*}{{\makecell[c]{}}} & {90 s}  &{0.5}  & {4.02/2C}  & {5.16\%}  & {35.237}\\
  \hline \hline
  \end{tabular}
\end{table}

\subsubsection{Optimal BESS Size with Different Parameter Settings}
\label{sec:SensitivityBESS}

\begin{figure}[tb]
  \centering
  \includegraphics[width=3.4in]{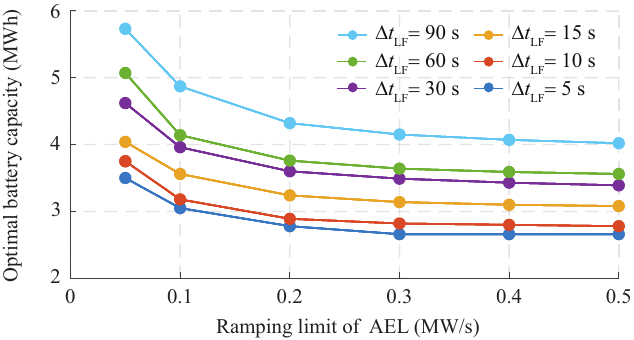}
  \vspace{-6pt}
  \caption{Optimal battery capacities with respect to different time steps of the LF and the ramping limit of the AEL.}
  \label{fig:SBatChangeRampRate}
\end{figure}

\begin{figure}[tb]
  \centering
  \includegraphics[width=3.5in]{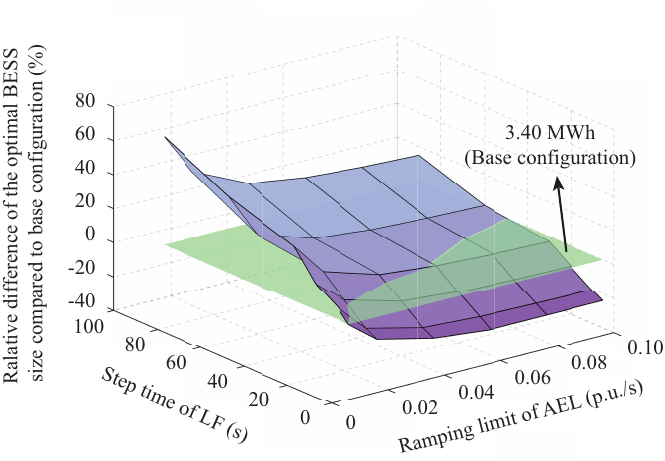}
  \vspace{-6pt}
  \caption{The ratio of the optimal battery capacities between the base parameter setting and other settings.}
  \label{fig:SBat3De}
\end{figure}

Figs. \ref{fig:SBatChangeRampRate} and \ref{fig:SBat3De} show the results of the sensitivity analysis. The results indicate that the optimal battery capacity increases almost exponentially when the ramping limit of AELs is lower than 0.3 MW/s (6\% of the rated power) across all configurations. When the ramping limit exceeds 0.3 MW/s, the battery capacity decreases slightly by 2\%, indicating that the ramping limit of 0.3 MW/s for AELs is a critical inflection point for battery capacity optimization.

Additionally, the reduction in the optimal battery capacity due to a decrease in the LF time step is more significant than that caused by an increase in the ramping limit. Numerically, as the time step and ramping limit increase to ten times their base values, the average reductions in the battery capacities are 0.876 MWh and 0.484 MWh, respectively. This demonstrates that the time step has a more substantial impact on sizing of the BESS. Therefore, with the proposed multi-timescale EMS, a shorter time step for the LF has a more pronounced effect than a higher load ramping limit of AELs. Meanwhile, under optimal sizing, the yearly degradation of the battery remains relatively constant at approximately 5\%, indicating that the BESS is fully utilized by the proposed EMS across different parameter settings.

An examination of Fig. \ref{fig:SBat3De} shows that the minimal battery capacity of 2.66 MWh is obtained when the ramping limit of AELs reaches 0.5 MW/s with a 5-s time step for the LF. However, it may be challenging to design the ramping at 0.5 MW/s (10\% of the rated power) for utility-scale AELs \cite{david2019advances}. Achieving such performance requires a high-bandwidth control loop in the EL rectifier, which may compromise control precision and complicate controller design. In addition, rapid load adjustments increase heat generation and hydrogen--oxygen gas crossover, necessitating larger heat exchangers or higher-power chillers and more advanced control for limiting hydrogen to oxygen (HTO) impurity. These enhancements incur additional costs and introduce design complexities. Considering current manufacturing capabilities of commercial ELs, a ramping limit of 0.2 MW/s to 0.3 MW/s (4\%--6\% of the rated power) with a 5-second or 10-second time step is more practical based on our experiments on several 5 MW-rated commercial AELs. Under these settings, the required BESS capacity remains within an acceptable range of 2.66 to 2.89 MWh.

\subsubsection{LCOH with Different Parameter Settings}
\label{sec:SensitivityLCOH}

\begin{figure}[tbp]
  \centering
  \includegraphics[width=3.5in]{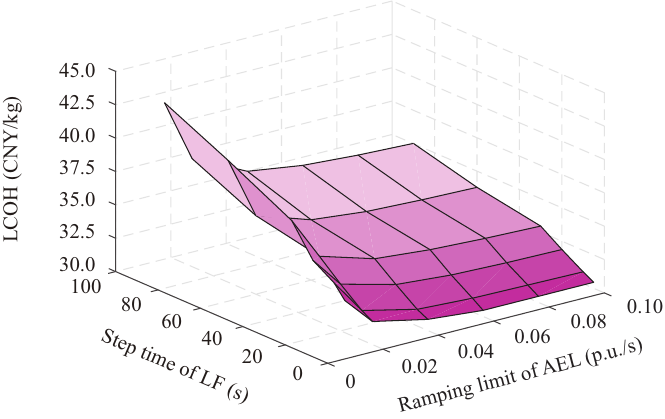}
  \vspace{-6pt}
  \caption{LCOH with respect to different time steps of the LF and ramping limit of the AEL.}
  \label{fig:LCOH}
\end{figure}

This section discusses the impact of various factors on the LCOH, with the results shown in Fig. \ref{fig:LCOH}.
It is observed that the cost of hydrogen decreases almost linearly with the decreasing time step of LF. Moreover, the rate of change of LCOH is significantly higher when the time step is shorter than 30 s. The ramping limit of AELs has a minimal impact on LCOH after it exceeds 0.3 MW/s. However, when the ramping limit drops below 0.3 MW/s, the LCOH rises exponentially, reaching its maximum value at 0.05 MW/s in all evaluated configurations.

Comparison with Fig. \ref{fig:LCOH} shows that the trend for the LCOH, which is influenced by both the ramping limit and the time step of LF, is similar to that of BESS size. This similarity indicates that the annual costs, including the capital cost and O\&M cost of BESS, significantly affect the LCOH. The minimal LCOH of 25.451 CNY/kg (3.510 USD/kg) appears when the ramping limit of AELs is 0.5 MW/s with a 5-s time step for the LF. When the ramping limit is set within the recommended range of 0.2 MW/s to 0.3 MW/s with a 5-s or 10-s time step for the LF, the LCOH also remains competitive, ranging from 25.458 CNY/kg (approx. 3.511 USD/kg) to 26.246 CNY/kg (approx. 3.620 USD/kg).

\subsection{Impact of Interannual Climate Variability on Optimal BESS Size and LCOH}
\label{sec:Conclusion}
\begin{figure}[tbp]
  \centering
  \includegraphics[width=3.5in]{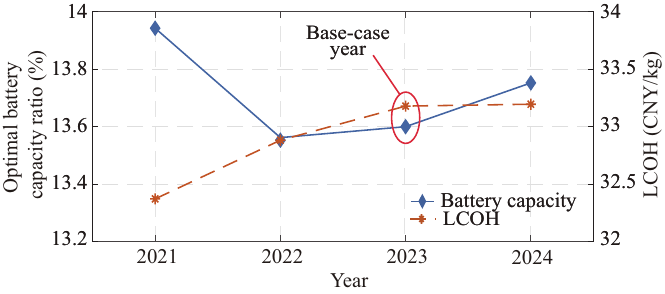}
  \vspace{-6pt}
  \caption{Optimal BESS size and LCOH with respect to different yearly meteorological data.}
  \label{fig:YearlyBESS}
\end{figure}

In order to analyze the influence of the inter-annual climate variability impact on the on BESS sizing results. This section re-evaluate the BESS size based on the base-case parameter settings using in Section \ref{sec:SimResult} with multi-year meteorological data. We define the \emph{optimal battery capacity ratio} as the percentage of battery capacity relative to the rated hourly energy output of renewable sources. As shown in Fig. \ref{fig:YearlyBESS}, although interannual climate variability has a certain influence on the optimal battery sizing, the impact remains limited. For example, due to slightly higher wind speeds in 2021, the optimal battery capacity ratio reached 13.96\%, whereas in other representative years, it remained consistently around 13.6\%.
Given the observed long-term declining trend of onshore wind speeds in China \cite{Tsinghua2024Outlook,CMA2023Bulletin}, a battery capacity ratio of 13.6\% is considered a reasonable and practical planning value. For more conservative system planning, configuring the battery capacity ratio at 14\% to accommodate extreme weather scenarios would be advisable.

\section{Conclusions and Outlook}
\label{sec:Conclusion}

\subsection{Conclusions}
To maintain reliable operation while ensuring the economic profitability of the OReP2HS, this work presents a multi-timescale EMS that coordinates the energy balance from transient power balance to intra-day scheduling for PV, WTs, grid-forming BESS, and AELs. The optimal BESS size under the proposed EMS is evaluated using a high-fidelity simulation-based search procedure, considering constraints related to the grid-forming capability, continuous operation during emergencies, and long-term energy balance that cover timescales ranging from milliseconds to a year.

A realistic OReP2HS planned in Inner Mongolia, China, is used as a case study. The studied system includes 18.75 MW WTs, a 6.25 MWp PV plant, and four AELs, each rated at 5 MW. The optimal size of the grid-forming BESS is evaluated to be 6.8 MW/3.4 MWh, with the battery capacity accounting for 13.6\% of the rated hourly energy output of the renewable power sources. However, the BESS exhibits a significant annual degradation rate of 4.87\%, increasing its O\&M costs for replacing retired batteries. Consequently, the capital expenditures of the BESS account for 17.83\% of the total, and the minimal LCOH under the base-case setting is evaluated to be 33.212 CNY/kg (4.581 USD/kg).

We also find that without considering voltage and frequency securities, the BESS size would be underestimated by approximately 15--30\% by the methods used in the studies reported in the literature, indicating that security constraints are necessary in the design of the OReP2HS. Additionally, the SOC of the BESS shows a positive correlation with renewable power. This implies that SOC can be predicted based on renewable power forecasts, enabling predictive measures to mitigate battery degradation.

Finally, sensitivity analysis shows that the required BESS size decreases when increasing the load ramping rate of AELs or shortening the time step of the LF control. Reducing the adjustment time step of the electrolytic load from 90 to 5 seconds and increasing its load ramping limit from 1\% to 10\% of the rated electrolytic load reduces the BESS size from 5.73 MWh to 2.66 MWh and the LCOH from 37.814 CNY/kg (approx. 5.216 USD/kg) to 25.458 CNY/kg (approx. 3.511 USD/kg). Considering the cost for the design, manufacture, and maintenance of utility-scale AELs with fast load regulation capability, an electrolytic load ramping limit at 4--6\% of the rated power and an adjustment time step at 5 or 10 seconds are recommended for the AEL to balance profitability and technological feasibility.

\subsection{\color{black}Methodological Generalizability }
{\color{black}The proposed multi-timescale EMS and the high-fidelity electromagnetic transient-based grid-forming BESS planning method are not limited to OReP2HSs. They can also be applied to a wide range of renewable and hydrogen-based energy systems, such as green buildings, remote communities, and off-grid industrial parks. Owing to their strong similarity in system configuration and operational characteristics to the studied OReP2HS, the proposed approach can be readily transferred and adapted to these scenarios.

Besides, the proposed BESS sizing framework is closely related to the multi-stage decision-making structure commonly adopted in process systems engineering \cite{brunaud2017perspectives,grossmann2019process}. Specifically, the EMS design that leverages the load flexibility of ELs to minimize the required BESS capacity corresponds to the first-stage, capital cost-oriented process design, while the rolling scheduling that aims to minimize operational costs represents the second-stage, operation-oriented decision-making. This alignment not only reinforces the methodological soundness of the proposed framework but also demonstrates its practical applicability and scalability for broader engineering implementations. }

\subsection{Outlook}
Even though this work proposes a multi-timescale EMS and investigates the optimal size of the grid-forming BESS required for an OReP2HS, several future research directions need to be explored.

First, as noted in Section~\ref{sec:EMS}, the proposed EMS still requires refinement. {\color{black}Future work will focus on incorporating the probabilistic characteristics of renewable power into rolling scheduling, enhancing the scheduling framework through stochastic or robust optimization to better handle renewable uncertainties, and integrating proactive BESS control actions to further improve system performance. Moreover, the EH strategy can be upgraded using advanced techniques such as predictive adaptive control or anomaly detection methods. Improving forecasting accuracy in load following may also further enhance the EMS performance.}

Second, the iterative searching procedure presented in Section \ref{sec:BatOptim} is designed for evaluation of the BESS size. However, for OReP2HS planning, the sizes of the ELs, power sources, and the electrical network also must be determined in some cases, making the dimension of the search space much larger. Given the complicated dynamics of the OReP2HS across timescales, a more efficient and mathematically advanced optimization approach is desired for future research.

Finally, industrial-scale OReP2HSs may be connected to downstream hydrogen buffer tanks, hydrogen delivery and refueling systems, and hydrogen consumers, such as green ammonia and methanol synthesis, melting, or refining. The proposed EMS, particularly for long-term energy balance and intra-day rolling scheduling, needs to be improved to cover these downstream sections; this is another important direction for future research.

\section*{Acknowledgement}

The authors gratefully acknowledge the financial support from the National Key Research and Development Program of China (2021YFB4000503) and the National Natural Science Foundation of China (52377116, 52377115, and 52307126).

\section*{Declaration of Interest}

None.

\section*{Data Availability}

The data related to this work are available upon request.

\appendix
\section{Details of the Simulation Models}
\label{sec:appA}

\setcounter{figure}{0}
\renewcommand{\thefigure}{A\arabic{figure}}
\setcounter{table}{0}
\renewcommand{\thetable}{A\arabic{table}}

The topologies and parameters of each component in both the refined and simplified simulation models are elaborated below.
\subsection{PV Plant}

\begin{figure}[H]
  \centering
  \includegraphics[width=5in]{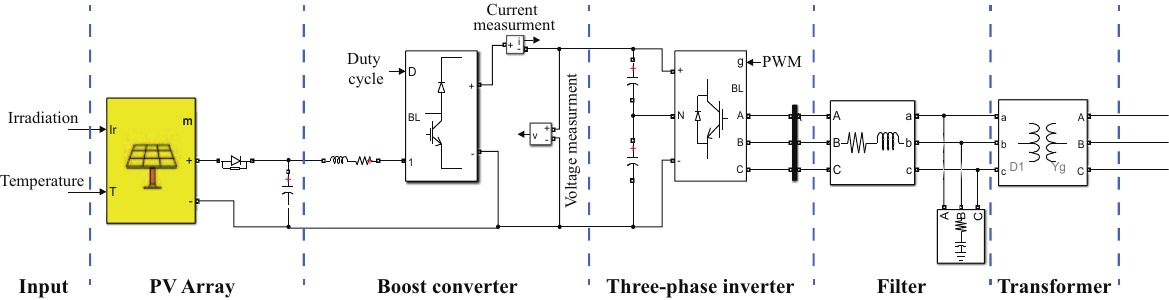}
  \vspace{-6pt}
  \caption{Topology of the PV plant model.}
  \label{fig:SolarPV}
\end{figure}

Fig. \ref{fig:SolarPV} depicts the PV plant model. It consists of a boost converter, a three-phase inverter, an AC filter, and a transformer connecting the PV array to the OReP2HS. The standard PV array module in the Simulink/Simscape environment is employed, and we set its parameters to meet the specifications of the PV plant in our studied case. The controllers and parameters of the PV plant model are shown in Table \ref{tab:ParaPV}.

\begin{table}[H]\footnotesize
  \renewcommand{\arraystretch}{1.0}
  \caption{Parameters of the PV plant model}
  \label{tab:ParaPV}
  \vspace{-6pt}
  \centering
  \begin{tabular}{ccc}
  \hline \hline
  Module  &  Parameter  & Value   \\ %
  \hline
  \multirow{7}{*}{{\makecell[c]{PV Array}}} & {Rated power of the PV array}  &{5.014 MW}  \\
  \multirow{7}{*}{{\makecell[c]{}}} & {Maximal power of per cell}  &{540.2 W}  \\
  \multirow{7}{*}{{\makecell[c]{}}} & {Open circuit voltage of per cell}  &{49.5 V}  \\
  \multirow{7}{*}{{\makecell[c]{}}} & {Short-circuit current of per cell}  &{13.85 A}  \\
  \multirow{7}{*}{{\makecell[c]{}}} & {Voltage at maximum power point}  &{41.65 V}  \\
  \multirow{7}{*}{{\makecell[c]{}}} & {Current at maximum power point}  &{12.97 A}  \\
  \multirow{7}{*}{{\makecell[c]{}}} & {Maximal power of per cell}  &{540.2 W}  \\
  \hline
  {Boost converter} & {Controller}  &{MPPT}  \\
  \hline
  \multirow{4}{*}{{\makecell[c]{Three-phase inverter}}} & {Rated DC Voltage}  &{1,500 V}  \\
  \multirow{4}{*}{{\makecell[c]{}}} & {Rated AC voltage}  &{630 V}  \\
  \multirow{4}{*}{{\makecell[c]{}}} & {Rated frequency}  &{50 Hz}  \\
  \multirow{4}{*}{{\makecell[c]{}}} & {Control}  &{DC voltage control}  \\
  \hline
  \multirow{3}{*}{{\makecell[c]{Transformer}}} & {Rated capacity and frequency}  &{5.0 MVA, 50 Hz}  \\
  \multirow{3}{*}{{\makecell[c]{}}} & {Primary voltage (RMS)}  &{35 kV}  \\
  \multirow{3}{*}{{\makecell[c]{}}} & {Secondary voltage (RMS)}  &{630 V}  \\
  \hline \hline
  \end{tabular}
\end{table}

%\break
\subsection{Wind Turbine}

\begin{figure}[H]
  \centering
  \includegraphics[width=4.5in]{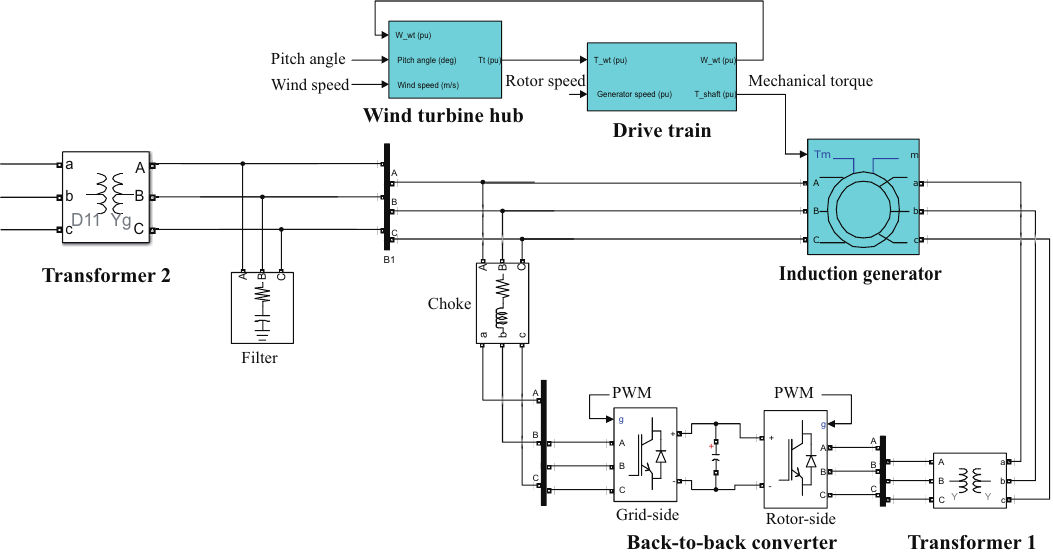}
  \vspace{-6pt}
  \caption{Topology of the wind turbine model.}
  \label{fig:WT}
\end{figure}

Fig. \ref{fig:WT} depicts the wind turbine model. A DFIG-based wind farm model in the Simulink/Simscape environment is employed, and its parameters are set to meet the specifications of the case study. The controllers and parameters of the wind turbine model are shown in Table \ref{tab:ParaWT}.

\begin{table}[H]\footnotesize
  \renewcommand{\arraystretch}{1.05}
  \caption{Parameters of the wind turbine model}
  \label{tab:ParaWT}
  \vspace{-6pt}
  \centering
  \begin{tabular}{ccc}
  \hline \hline
  Module  &  Parameter  & Value   \\ %
  \hline
  \multirow{6}{*}{{\makecell[c]{Wind turbine hub}}} & {Rated power}  &{6.25 MW}  \\
  \multirow{6}{*}{{\makecell[c]{}}} & {Impeller diameter}  &{175 m}  \\
  \multirow{6}{*}{{\makecell[c]{}}} & {Swept area}  &{24,053 m$^{3}$}  \\
  \multirow{6}{*}{{\makecell[c]{}}} & {Rated speed}  &{10.2 m/s}  \\
  \multirow{6}{*}{{\makecell[c]{}}} & {Cut-in speed}  &{3 m/s}  \\
  \multirow{6}{*}{{\makecell[c]{}}} & {Cut-off speed}  &{24 m/s}  \\

  \hline
  \multirow{2}{*}{{\makecell[c]{Drive train}}} & {Inertia constant}  &{5.25 s}  \\
  \multirow{2}{*}{{\makecell[c]{}}} & {model}  &{2-mass model}  \\

  \hline
  \multirow{3}{*}{{\makecell[c]{Induction Generator}}} & {Rated capacity}  &{6.25 MVA}  \\
  \multirow{3}{*}{{\makecell[c]{}}} & {Rated voltage}  &{690 V}  \\
  \multirow{3}{*}{{\makecell[c]{}}} & {Rated frequency}  &{50 Hz}  \\

  \hline
  \multirow{5}{*}{{\makecell[c]{Back-to-back converter}}} & {Rated capacity of the grid side}  &{1084 kVA}  \\
  \multirow{5}{*}{{\makecell[c]{}}} & {Rated capacity of the rotor side}  &{2,174 kVA}  \\
  \multirow{5}{*}{{\makecell[c]{}}} & {Rated AC voltage}  &{1,140 V}  \\
  \multirow{5}{*}{{\makecell[c]{}}} & {Grid-side controller}  &{DC voltage control}  \\
  \multirow{5}{*}{{\makecell[c]{}}} & {Rotor-side controller}  &{P/Q control}  \\

  \hline
  \multirow{3}{*}{{\makecell[c]{Transformer 1}}} & {Rated capacity and frequency}  &{6.25 MVA, 50 Hz}  \\
  \multirow{3}{*}{{\makecell[c]{}}} & {Primary voltage (RMS)}  &{1970 V}  \\
  \multirow{3}{*}{{\makecell[c]{}}} & {Secondary voltage (RMS)}  &{690 V}  \\

  \hline
  \multirow{3}{*}{{\makecell[c]{Transformer 2}}} & {Rated capacity and frequency}  &{6.25 MVA, 50 Hz}  \\
  \multirow{3}{*}{{\makecell[c]{}}} & {Primary voltage (RMS)}  &{35 kV}  \\
  \multirow{3}{*}{{\makecell[c]{}}} & {Secondary voltage (RMS)}  &{1140 V}  \\

  \hline \hline

  \end{tabular}
\end{table}

\subsection{BESS}

\begin{figure}[H]
  \centering
  \includegraphics[width=5in]{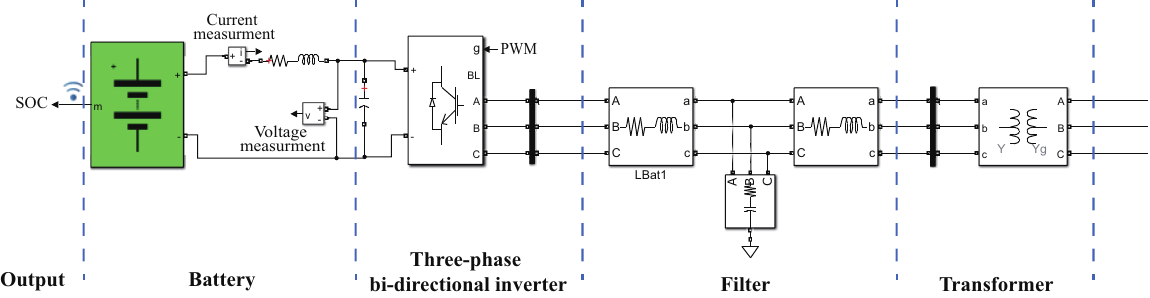}
  \caption{Topology of the BESS model.}
  \label{fig:Battery}
\end{figure}

The BESS model is depicted in Fig. \ref{fig:Battery}. It includes  a three-phase bi-directional inverter, an AC filter, and a transformer. The battery itself is modeled as a controllable ideal DC source in series with an internal resistance. The controllers and parameters of the BESS are summarized in Table \ref{tab:ParaBESS}.

\begin{table}[H]\footnotesize
  \renewcommand{\arraystretch}{1.05}
  \caption{Parameters of the BESS model}
  \label{tab:ParaBESS}
  \vspace{-6pt}
  \centering
  \begin{tabular}{ccc}
  \hline \hline
  Module  &  Parameter  & Value   \\ %
  \hline
  \multirow{6}{*}{{\makecell[c]{Battery}}} & {Rated capacity of the battery}  &{3.4 MWh}  \\
  \multirow{6}{*}{{\makecell[c]{}}} & {Charging/discharging rate}  &{2C}  \\
  \multirow{6}{*}{{\makecell[c]{}}} & {Rated capacity of per cell}  &{280 Ah}  \\
  \multirow{6}{*}{{\makecell[c]{}}} & {Rated Voltage of per cell}  &{3.2 V}  \\
  \multirow{6}{*}{{\makecell[c]{}}} & {Operating voltage of per cell}  &{2.8 V$\sim$4.0 V}  \\
  \multirow{6}{*}{{\makecell[c]{}}} & {Internal resistance}  &{0.39 M$\Omega$}  \\

  \hline
  \multirow{4}{*}{{\makecell[c]{Three-phase\\bi-directional inverter}}} & {Rated DC Voltage}  &{1300 V}  \\
  \multirow{4}{*}{{\makecell[c]{}}} & {Rated AC voltage}  &{690 V}  \\
  \multirow{4}{*}{{\makecell[c]{}}} & {Rated frequency}  &{50 Hz}  \\
  \multirow{4}{*}{{\makecell[c]{}}} & {Controller}  &{V-F droop control}  \\

  \hline
  \multirow{3}{*}{{\makecell[c]{Transformer}}} & {Rated capacity and frequency}  &{3.5 MVA, 50 Hz}  \\
  \multirow{3}{*}{{\makecell[c]{}}} & {Primary voltage (RMS)}  &{35 kV}  \\
  \multirow{3}{*}{{\makecell[c]{}}} & {Secondary voltage (RMS)}  &{690 V}  \\

  \hline \hline

  \end{tabular}
\end{table}

%\break
\subsection{AEL}

\begin{figure}[H]
  \centering
  \includegraphics[width=5in]{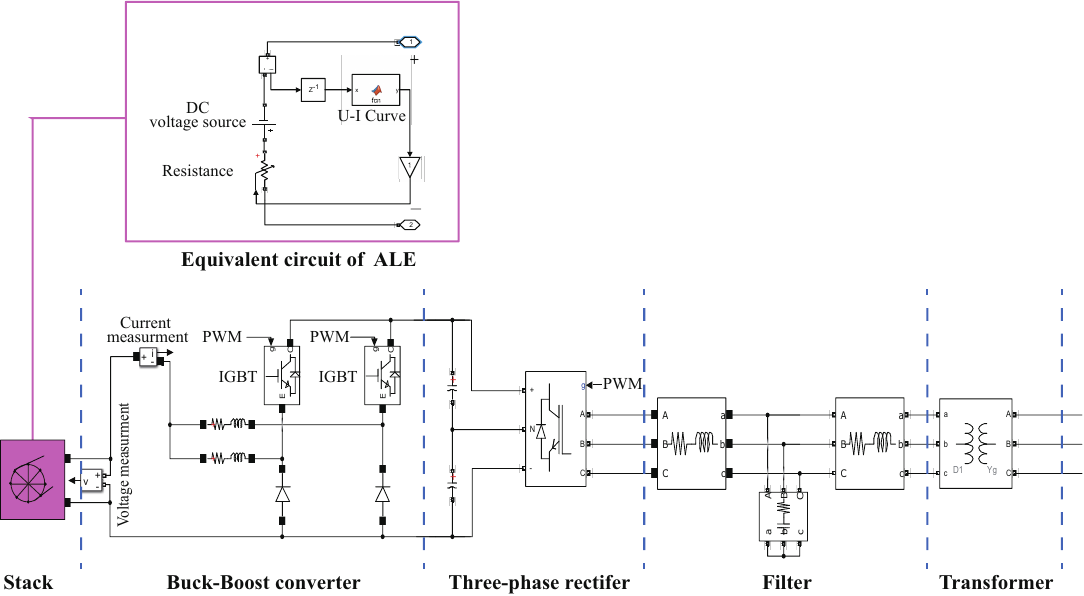}
  \caption{Topology of the AEL model.}
  \label{fig:ALE}
\end{figure}

\begin{figure}[htbp]
  \centering
  \includegraphics[width=3.0in]{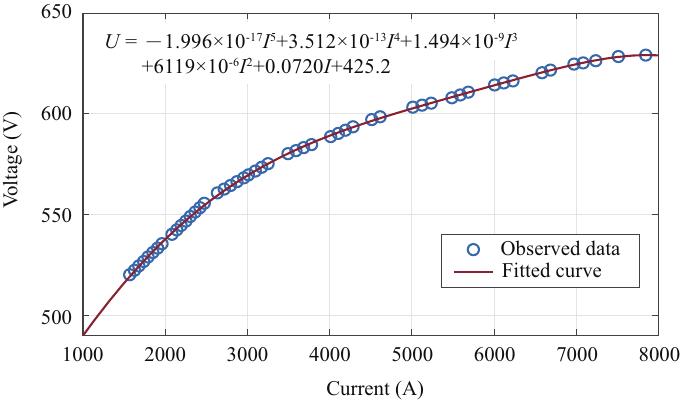}
  \vspace{-6pt}
  \caption{The experimental and fitted U-I curve of a 1,000 Nm$^3/$h-rated electrolysis stack produced by Peric Hydrogen Technologies.}
  \label{fig:UICurve}
\end{figure}

The AEL model is depicted in Fig. \ref{fig:ALE}. It includes a stack, a Buck/Boost converter, a three-phase rectifier, an AC filter, and a transformer connecting the AEL to the OReP2HS. The stack is represented by an equivalent circuit consisting of an nonlinear internal resistance in series with a DC voltage source to simulate its U-I characteristics, as shown in Fig. \ref{fig:UICurve}. {\color{black}The impact of the EDL on stack load adjustment is emulated by introducing a load ramping limit mechanism in this study.} The controllers and parameters of the AEL model are summarized in Table \ref{tab:ParaALE}.

\begin{table}[H]\footnotesize
  \renewcommand{\arraystretch}{1.}
  \caption{Parameters of the AEL model}
  \label{tab:ParaALE}
  \vspace{-6pt}
  \centering
  \begin{tabular}{ccc}
  \hline \hline
  Module  &  Parameter  & Value   \\ %
  \hline
  \multirow{2}{*}{{\makecell[c]{Equivalent circuit of AEL}}} & {Resistance}  &{0.05 $\Omega$}  \\
  \multirow{2}{*}{{\makecell[c]{}}} & {U-I curve}  &{Depicted in Fig. \ref{fig:UICurve}}  \\

  \hline
  {Buck-Boost converter} & {Controller}  &{DC current control}  \\

  \hline
  \multirow{4}{*}{{\makecell[c]{Three-phase inverter}}} & {Rated capacity }  &{5.0 MW}  \\
  \multirow{4}{*}{{\makecell[c]{}}} & {Rated AC voltage}  &{750 V}  \\
  \multirow{4}{*}{{\makecell[c]{}}} & {Controller}  &{DC voltage control}  \\

  \hline
  \multirow{3}{*}{{\makecell[c]{Transformer}}} & {Rated capacity and frequency}  &{6.0 MVA, 50 Hz}  \\
  \multirow{3}{*}{{\makecell[c]{}}} & {Primary voltage (RMS)}  &{35 kV}  \\
  \multirow{3}{*}{{\makecell[c]{}}} & {Secondary voltage (RMS)}  &{380 V}  \\

  \hline \hline

  \end{tabular}
\end{table}

%\break

\subsection{\color{black}Selection of the Time-Step For Simulation}
{\color{black}As the simulation time step influence both simulation accuracy and computational efficiency. This section provide a time-step sensitivity analysis based on the proposed high-fidelity model over a 10-hour production simulation with three different time-step settings, $\Delta t \in \{0.02, 0.04, 0.08\}$~ms. The comparison of key converter and system indicators, including the DC-link capacitor voltage of the AC/DC converter of BESS ($v_{dc}$) , effective value of AC current of BESS ($i_{\mathrm{AC,RMS}}$), effective value of  AC voltage of BESS ($v_{\mathrm{AC,RMS}}$), and system frequency $f$, is summarized in Table~\ref{tab:StepSensitivity}.
\begin{table}[h]\footnotesize
  \renewcommand{\arraystretch}{1.0}
  \caption{\color{black}Sensitivity analysis of the of the time step for simulation.}
  \label{tab:StepSensitivity}
  \vspace{-6pt}
  \centering
  {\color{black}
  \setlength{\tabcolsep}{3.2pt} % 控制列间距
  \begin{tabular}{cccccccccc}
    \hline\hline
    \multirow{2}{*}{$\Delta t$ (ms)} & \multirow{2}{*}{Computation time (h)}
      & \multicolumn{2}{c}{$v_{dc}$}
      & \multicolumn{2}{c}{$i_{\mathrm{AC,RMS}}$}
      & \multicolumn{2}{c}{$v_{\mathrm{AC,RMS}}$}
      & \multicolumn{2}{c}{$f$} \\
    \cline{3-10}
      &  & MSE & MAX & MSE & MAX & MSE & MAX & MSE & MAX \\
    \hline
     0.02 (ref) &  36.42 & / & / & / & / & / & / & / & / \\
    0.04       &   20.12 & $1.5\times10^{-4}$ & $0.07\%$ & $4.4\times10^{-4}$ & $0.13\%$ & $5.8\times10^{-4}$ & $0.17\%$ & $6.9\times10^{-5}$ & $0.009$\,Hz \\
    0.08       &   16.94 & $3.3\times10^{-3}$ & $0.41\%$ & $3.7\times10^{-3}$ & $0.65\%$ & $4.5\times10^{-3}$ & $0.84\%$ & $1.4\times10^{-4}$ & $0.017$\,Hz \\
    \hline\hline
  \end{tabular}
  \vspace{2pt}}
  \begin{minipage}{0.95\linewidth}\footnotesize
  \color{black}Notes: MSE and MAX are calculated with respect to the 0.02 ms reference solution. RMS quantities are computed using a one-cycle sliding window at 50\,Hz. The frequency $f$ is estimated from the PLL output. All candidate waveforms are time-aligned and interpolated to the reference time grid before comparison.
  \end{minipage}
\end{table}
The mean-square error (MSE) and maximum deviation (MAX) used in Table~\ref{tab:StepSensitivity} are defined as
\begin{equation}
\mathrm{MSE}(x) = \frac{1}{N}\sum_{k=1}^{N}\!\big(x(k)-x_{\mathrm{ref}}(k)\big)^{2},
\qquad
\mathrm{MAX}(x) = \max_{k}\!\big|x(k)-x_{\mathrm{ref}}(k)\big|,
\end{equation}
where $x_{\mathrm{ref}}(k)$ is the reference waveform obtained at $\Delta t_{\mathrm{ref}}=0.02\,\mathrm{ms}$, and $x(k)$ is the waveform computed at a candidate time step $\Delta t$.

The results indicate that $\Delta t=0.04$ ms reproduces the reference results with sub-percent deviations in all indicators, while reducing computation time by nearly 44.76\%. Although using $\Delta t=0.08$ ms further decreases the computation time slightly, the model accuracy deteriorates by approximately one order of magnitude. Therefore, 0.04 ms shows a balance that ensures high-fidelity simulation of converter switching dynamics without incurring excessive computational cost in the simulations, thus be selected as the time step in our study.}

%\break

\section{Nomenclature}
\label{sec:appB}

\setcounter{figure}{0}
\renewcommand{\thefigure}{B\arabic{figure}}
\setcounter{table}{0}
\renewcommand{\thetable}{B\arabic{table}}

%\begin{longtable}[htbp]%\footnotesize
%  \renewcommand{\arraystretch}{1}
%  \centering
  \begin{longtable}{p{1.8cm}p{5.6cm}p{1.8cm}p{5.6cm}}
  \emph{Abbreviations}     &     &     &  \\
  AC      &   Alternating current  & $\Delta t_{\text{RS}}$           & Time step of the scheduling \\
  EL     &  Electrolyzer  & $\Delta t_\text{{LF}}$           & Time step of LF control \\
  BESS     &  Battery energy storage system  & $\Delta E_{j}^{\mathrm{degra,max}}$           & Battery degradation limit \\
  %CODE     &  Continuous operation during emergencies  & $\alpha$           & Weigh of the MA. \\
  EMS     &  Energy management system  & $\beta$           & Intercept of SOC correction control \\
  LCOH    &  Levelized cost of hydrogen  & $\lambda_{j}^{\mathrm{O\&M}}$           & Ratio of annual O$\&$M cost\\
  MA      &  Moving average  & $\eta^{\mathrm{BES}}$           & BESS Charge/discharge efficiency \\
  MILP      &  Mixed-integer linear programming  & \emph{Variables}        &  \\
  MPPT      &  Maximum power point tracking  &   $C_{\mathrm{init}}$    &  Initial investment cost \\
  OReP2HS      &  Off-grid renewable P2H system  &  $C^{{\mathrm{fixed}}}_{\mathrm{O\&M}} $     & Fixed annual O$\&$M cost \\
  PCC      &  Point of common connection  &  $C^{\mathrm{rep}}_{j} $         &  Battery replacement cost \\
  PV      & Photovoltaic &  $C^{\mathrm{rec}}_{j} $         &  Battery recycling revenue \\
  ReP2H     & Renewable power to hydrogen &  $C^{{\mathrm{vari}}}_{\mathrm{O\&M}} $         &  Variable annual O$\&$M cost \\
  RoCoF     & Rate of change of frequency &  $iter $         &  Iterator in battery size search \\
  O$\&$M     & Operation and maintenance & $M^{\mathrm{H_{2}}} $         & Annual hydrogen yield \\
  SOC     & State of charge of battery & $K^{\mathrm{rep}}_{j} $         & Battery replacement times \\
  LF     & Load following & $P^{\mathrm{EL}}_{i,t} $         & Load power of AEL \\
  WT     & Wind turbine  & $P^{\mathrm{EL}}_{\mathrm{sb}} $         & Standby power of EL \\
   \emph{Parameters}     &   &  $P^{\mathrm{EL,RS}}_{i,t} $       & Baseline load command of EL \\
  $C^{\mathrm{EL}}_{\mathrm{up,h/c}}$    &  Hot/cold startup costs of EL   &  $P^{\mathrm{EL,LF}}_{i,t} $       & Load correction command of EL \\
  $C^{\mathrm{EL}}_{\mathrm{down}}$    &  Costs for hot/cold start-up of EL  &  $P^{\mathrm{BES,C/D}}_{t}$       & BESS Charging/discharging power  \\
   $C^{\mathrm{H_{2}}}$    &  Selling price of hydrogen  &  $\tilde{P}_\tau^\mathrm{RES}$      & Fast prediction in LF control \\
  $C^{\mathrm{rep}}_{j} $  &  Unit replacement cost of device $j$  &   $  P^{\mathrm{WT/PV,cut}}_{i,t}$  & Power curtailment of WT/PV \\
  $C^{\mathrm{rec}}_{j} $  &  Unit recycling revenue of device $j$  &  $\widehat P^{\mathrm{WT}/\mathrm{PV}}_{i,t} $    & Ultra-short-term forecast of WT/PV  \\
  $C^{\mathrm{unit}}_{j} $  &  The unit cost of device $j$  &  $ S_{j} $    & Capacity of device $j$  \\
  $f^{\mathrm{PCC}}_{\text{min/min}} $  & Frequency limits   &  $ S^{\mathrm{rep}}_{j}$    & Replacement capacity of device $j$   \\
  $K^{\mathrm{H_{2}}}_{i} $  &  Energy conversion coefficient of AEL  &  $ S^{\mathrm{rec}}_{j}$    & Recycling capacity of device $j$   \\
  $ r $  &  Discount rate  &  $ SOC^{\mathrm{BES}}_{t}$    & State of charge of battery   \\
  $ S^{\mathrm{BES}}_{\mathrm{init}}$   & Initial size of BESS in the search  &  $ S^{\mathrm{EL}}_{i}$    & Capacity of electrolyzer   \\
  $ SOC^{\mathrm{BES}}_{\text{min}/\text{max}}$   &  Battery state of charge limits &  $ S^{\mathrm{ds},\text{iter}}_{\mathrm{BES}}$    & Iterator in battery sizing   \\
  $ s  $   & Penalty for power curtailment  &  $ T^{\mathrm{Fail}}_{j}$    & Durable years of device $j$   \\
  $ T_{\mathrm{RS}}$   & Total time periods of scheduling  &  $ \mu^{\mathrm{BES,C/D}}_{t}$    & BESS Charging/discharging state    \\
  $ T^{\mathrm{down}}_{\text{min}}$   & Minimal restart interval for AEL  &  $ \mu^{\mathrm{st}}_{i,t}/\mu^{\mathrm{sb}}_{i,t}/\mu^{\mathrm{sp}}_{i,t} $    & Production/standby/shut-down state  of AEL \\
  $ T^{\mathrm{LCC}}_{j}$   & Lifetime of device $j$  &  $ \mu^{\mathrm{up,h}}_{i,t}/\mu^{\mathrm{up,c}}_{i,t}$/ $\mu^{\mathrm{down}}_{i,t} $    & Action of hot start-up/cold start-up/shut-down of AEL  \\
  $V^{\mathrm{PCC}}_{\text{min/max}}$  &  Voltage limits   &  $ \Delta E^{\mathrm{degra,year}}_{j}$    & Annual degradation of BESS   \\
  $r^{\mathrm{EL}}_{\text{min/max}} $  & Load regulation limits of EL   &  &  \\%\emph{ Subscripts } &  \\

  \end{longtable}
%\end{longtable}
%\break
%\section*{References}

%\bibliographystyle{elsarticle-num}
%\bibliography{EMS_BatOptim}

\end{document}